\documentclass[preprint,12pt]{elsarticle}

\pdfoutput=1
\usepackage{amssymb}

\usepackage{epsfig}
\usepackage{amsmath}
\usepackage{amssymb}
\usepackage{amsthm}
\usepackage{graphicx}
\usepackage{caption}
\usepackage{epic}
\usepackage{latexsym}
\usepackage[usenames]{color}
\newtheorem{theorem}{Theorem}
\newtheorem{definition}{Definition}
\newtheorem{lemma}{Lemma}
\newtheorem{corollary}{Corollary}

\numberwithin{theorem}{section}
\numberwithin{definition}{section}
\numberwithin{corollary}{section}
\numberwithin{lemma}{section}
\numberwithin{equation}{section}
\renewcommand{\thesection}{\arabic{section}}

\renewcommand{\figurename}{Fig.}
\renewcommand{\tablename}{Table}
\everymath{\displaystyle}
\def\l{\left}
\def\r{\right}

\def\R{\mathbb{R}}
\def\Z{\mathbb{Z}}
\def\O{\mathcal{O}}
\def\qed{\hspace{3mm}\hbox{\vrule height7pt width5pt}}
\newcommand{\bb}[1]{\begin{equation}\label{#1}}
\newcommand{\ee}{\end{equation}}
\newcommand{\bbb}{\begin{eqnarray}}
\newcommand{\eee}{\end{eqnarray}}
\newcommand{\bbbb}{\begin{eqnarray*}}
\newcommand{\eeee}{\end{eqnarray*}}
\newcommand{\nnn}{\nonumber}
\newcommand{\tT}{\intercal}
\newcommand{\red}{\color{red}}

\newcommand{\black}{\color{black}}

\journal{Journal of Computational and Applied Mathematics}
\begin{document}

\begin{frontmatter}

\title{A semi-adaptive finite difference method for simulating two-sided fractional convection-diffusion quenching problems }

\author[1]{Rumin Dong\fnref{fn1}} \ead{dongrumin0113@163.com}
\author[1]{Lin Zhu\corref{cor1}\fnref{fn1}} \ead{zhulin\_nx@163.com}
\author[2]{Qin Sheng\fnref{fn2}} \ead{qin\_sheng@baylor.edu}
\author[1]{Bingxin Zhao\fnref{fn1}} \ead{zhao\_bx@nxu.edu.cn}

\cortext[cor1]{Corresponding author} \fntext[fn1]{This author is supported in part by the
National Science Foundation of China (No. 12062024)
and Funds for Ningxia University Scientific Research Projects (Grant No. NYG2024039), China.}
\fntext[fn2]{This author is supported in part by the National Science Foundation (Grant No. DMS-2318032) and
Simons Foundation (Grant No. MPS-1001466), USA.}
\address[1]{School of Mathematics and Statistics, Ningxia Key Laboratory of Interdisciplinary Mechanics and
Scientific Computing, Ningxia Basic Science Research Center of Mathematics, Ningxia University, Yinchuan 750021, China}
\address[2]{Department of Mathematics and Center for Astrophysics, Space Physics and Engineering
Research, Baylor University, Waco, TX 76798-7328, USA}

\begin{abstract}
This paper investigates quenching solutions of an one-dimensional, two-sided Riemann-Liouville fractional order
convection-diffusion problem. Fractional order spatial derivatives are discretized using weighted averaging
approximations in conjunction with standard and shifted Gr\"{u}nwald formulas. The advective term is handled
utilizing a straightforward Euler formula, resulting in a semi-discretized system of nonlinear ordinary differential
equations. The conservativeness of the proposed scheme is rigorously proved and validated through simulation
experiments. The study is further advanced to a fully discretized, semi-adaptive finite difference method.
Detailed analysis is implemented for the monotonicity, positivity and stability of the scheme. Investigations are carried
out to assess the potential impacts of the fractional order on quenching location, quenching time, and critical length.
The computational results are thoroughly discussed and analyzed, providing a more comprehensive
understanding of the quenching phenomena modeled through two-sided fractional order convection-diffusion problems.

\end{abstract}




\end{frontmatter}




\section{Introduction}\label{sec:  sec1}

Fractional order derivatives have been playing a significant role in various fields due to their ability to model certain complex
behaviors more accurately than traditional integer order derivatives. The nonlocal nature of fractional derivatives allows
for the capture of memory and hereditary properties in materials and systems, making them particularly useful in
scenarios such as viscoelastic material designs \cite{chen2014flutter}, digital signal processing \cite{calderon2006fractional},
optimal controls \cite{li2014synchronization}, biomedicine \cite{anastasio1994fractional}, and
water resource preservation \cite{benson2000application,benson2000fractional}. Utilizing the nonlocal capabilities,
fractional partial differential equations (FPDEs) have become increasingly effective tools to analyze
systems across many diverse applications \cite{oldham1974fractional}. Various numerical strategies have been
developed to handle the unique challenges posed by FPDEs. These advances are crucial for leveraging the
full potential of fractional calculus in modeling and solving real-world problems.



Among the numerous numerical strategies implemented for solving FPDEs in recent decades, adaptive finite
difference methods \cite{zhang2009finite,zhu2020note,sheng2023nonconventional,sk2001,padgett2018quenching},
finite element methods \cite{agrawal2008general}, spectral methods \cite{lin2007finite}, finite volume methods \cite{feng2015stability},
and differential transformation methods \cite{biazar2011differential} are particularly effective. Among them,
a special finite volume strategy for noninteger order diffusion equations was developed by Liu et al. \cite{liu2014new}.
Furthermore, Hejazi et al. proposed an effective methodology for assessing the stability and convergency of numerical
solutions from fractional order problems \cite{hejazi2014stability}.
Rui and Huang later constructed a stable finite difference scheme for certain FPDEs \cite{rui2015uniformly}. On the other
hand,  there have been numerous discussions in applications of the meshless and machine learning methods for
approximating solutions of FPDEs. For instance, Badr et al. demonstrated an optimized spline coloration strategy for solving
temporal FPDEs and its stability was investigated systematically \cite{badr2018stability}.
Zafarghandi et al. introduced radial basis function equipped collocation techniques for Riesz space fractional
advection-dispersion equations and detailed analysis was provided in \cite{zafarghandi2019numerical}. Further,
Du et al. established an innovative meshless method to calculate minimum residual approximations of the solution of a multi-term,
time-fractional, integro-differential equation \cite{du2023meshless}.

While quenching phenomena have been well observed in fields of galaxy evolution, laser propagations,
nanofluids and superconductivities \cite{babu2017quenching,wang2021estimation,martinez2022nonlinear,ren2023study},
numerical analysis and simulations of singular quenching solutions of modeling equations are still in infancy.
In this highly valuable new territory, Sheng and Khaliq proceeded and
comprised an adaptive algorithm for nonlinear fractional order problems and modeled critical limits
corresponding to quenching phenomena \cite{sheng1999compound}.
Padgett further investigated a time-space fractional semi-linear quenching model.
Necessary conditions for quenching to occur were established \cite{padgett2018quenching}.
To ensure an overall numerical stability and take care of
the quenching nonlinearity, Sheng and Torres published detailed stability analysis without freezing the nonlinear source
entries of principal quenching differential equations \cite{sheng2023nonconventional,torres2024}. Following the ideas,
Zhu et al. studies quenching phenomena observed in fractional order diffusions modeled through one-sided and two-sided
Riemann-Liouville derivatives, respectively \cite{zhu2020note,zhu2023simulation}.
In 2023, Liu, Zhu and Sheng focused and analyzed a series of computational investigations
of fractional order convection-diffusion FPDEs in one-dimensional, one-sided fashion \cite{liu2023semi}. Interesting
results and discoveries acquired inspired the current project.

This paper aims at quenching phenomena occurred in one-dimensional, two-sided non-integer
order convection-diffusion equations.
The study consists of following parts.
In Section 2, we propose a semi-discrete adaptive scheme for a bilateral fractional-order convection-diffusion
problem with nonlinear source term. In Section 3, analysis of the monotonicity and positivity of the finite
difference scheme are conducted. The numerical method is proved to be stable.
Section 4 is devoted to simulations of anticipated critical lengths, quenching time and quenching locations
via four interconnected numerical experiments. A fully discretized semi-adaptive method is built, analyzed and
used. The convergence is numerically validated via simulations.
Section 5 summarizes observations and analysis about the numerical procedures accomplished.
Continuing research efforts and directions are reported.

\section{Semi-adaptive fractional model}
\label{sec: sec2}

This paper concerns the semi-linear nonlocal convection-diffusion Kawarada problem:
\bbb
\frac{\partial v(x,t)}{\partial t}&=&d_{+}(x,t)\frac{\partial^{\sigma}v(x,t)}{\partial_{+}x^{\sigma}}+d_{-}(x,t)\frac{\partial^{\sigma}v(x,t)}
{\partial_{-}x^{\sigma}}+c(x,t)\frac{\partial v(x,t)}{\partial x}+p(v), \nnn\\
v(0,t)&=&v(a,t)~=~0, \label{e2-1}\\
v(x,0)&=&\psi(x),\nnn
\eee
where
$\frac{\partial^{\sigma}v(x,t)}{\partial_{+}x^{\sigma}}$ and $\frac{\partial^{\sigma}v(x,t)}{\partial_{-}x^{\sigma}},~1<\sigma\leq 2,$
are $\sigma$-th order spatial partial derivatives of the right and left in the Riemann-Liouville sense
\cite{oldham1974fractional,zhu2020note,zhu2023simulation}, respectively,
$d_{+}(x,t)>0,$ $d_{-}(x,t)>0,$ $\kappa> \psi\geq 0,~0\leq x\leq a,\, 0\leq t<T,$ and $p(v)=(\kappa-v)^{-\theta},~\theta>0.$

Set $h=\frac{a}{L}, \tau=\frac{T}{M}$ for $L,M\in \Z^{+}.$ We define a space-time mesh region $\mathbb{\R}_{L,M} = \l\{(x_k,t_n)\r\}_{L,M},$ where
$$
x_{k}=kh,~k=0,1,2,\ldots,L;~
t_{n}=n\tau,~n=0,1,2,\ldots,M.
$$
To achieve a second-order accuracy in the space, we employ following standard and weighted average of shifted
Gr\"{u}nwald formulas for approximating the spatial derivative quantities
$\frac{\partial^{\sigma}v(x_k,t)}{\partial_{+}x^{\sigma}}$ and $\frac{\partial^{\sigma}v(x_k,t)}{\partial_{-}x^{\sigma}},$
respectively:
\bbb
\frac{\partial^{\sigma}v(x_{k},t)}{\partial_{+}x^{\sigma}}&=&\frac{1}{h^{\sigma}}\l[\l(1-\frac{\sigma}{2}\r)\sum_{j=0}^{k}z_{j}v_{k-j}+
\frac{\sigma}{2}\sum_{j=0}^{k+1}z_{j}v_{k-j+1}\r]+\O(h^{2});\label{e2-2}\\
\frac{\partial^{\sigma}v(x_{k},t)}{\partial_{-}x^{\sigma}}&=&\frac{1}{h^{\sigma}}\l[\l(1-\frac{\sigma}{2}\r)\sum_{j=0}^{L-k}z_{j}v_{k+j}+
\frac{\sigma}{2}\sum_{j=0}^{L-k+1}z_{j}v_{k+j-1}\r]+\O(h^{2}),\label{e2-3}
\eee
where $v_{k+j} =v(x_{k+j},t)$ as $h\rightarrow 0^+,$ and
\bb{e2-4}
\left\{
  \begin{array}{ll}
    z_{0}=1,~z_{1}=-\sigma<0,~z_{2}=\frac{\sigma(\sigma-1)}{2}>0,\\
    1\geq z_{2}\geq z_{3}\geq\cdots\geq0, \\
    \sum_{m=0}^{\infty}z_{m}=0,~\sum_{m=0}^{K}z_{m}\leq0,~~~K\in \Z^{+}.
 \end{array}
\right.
\ee

We further have
\bb{e2-5}
\frac{\partial v(x_{k},t)}{\partial x} 
=\frac{v_{k+1}-v_{k}}{h}+\O(h),~~~h\rightarrow 0^{+}.
\ee

Thus, our semi-discretized differential system can be obtained by applying (\ref{e2-2})-(\ref{e2-5}) to (\ref{e2-1}): 
\bbb
(v_{t})_{k}
&=&\frac{d_{+,k}}{h^{\sigma}}\l[\l(1-\frac{\sigma}{2}\r)\sum_{j=0}^{k}z_{j}v_{k-j}+\frac{\sigma}{2}\sum_{j=0}^{k+1}z_{j}v_{k-j+1}\r]
+\frac{d_{-,k}}{h^{\sigma}}\l[\l(1-\frac{\sigma}{2}\r)\sum_{j=0}^{L-k}z_{j}v_{k+j}\r.\nnn\\
&&\l.+\frac{\sigma}{2}\sum_{j=0}^{L-k+1}z_{j}v_{k+j-1}\r]+\frac{c_{k}}{h}(v_{k+1}-v_{k})+p_{k},~~k=0,1,2,\ldots,L,
\eee
where $d_{+,k}=d_{+}(x_{k},t),\; d_{-,k}=d_{-}(x_{k},t)$.
Combining the boundary conditions in (\ref{e2-1}), we formulate our semidiscretized scheme to
\bb{e2-7}
\frac{dv}{dt}=Sv+p(v),~0<t<T;~v(0)=v_{0},
\ee
where
$v=(v_{1},v_{2},\cdots,v_{L})^{\tT}\in\R^{L},~
p=(p_{1},p_{2},\cdots,p_{L})^{\tT}\in\R^{L},~
S=\l[s_{k,j}\r]\in\R^{L\times L}.$
Further,
\bbbb
s_{k,j}=\left\{
         \begin{array}{ll}
           \frac{d_{+,k}+d_{-,k}}{h^{\sigma}}\l[\l(1-\frac{\sigma}{2}\r)z_{0}+\frac{\sigma}{2}z_{1}\r]-\frac{c_{k}}{h},
& k=j, \hbox{ } \\
           \frac{d_{+,k}\sigma}{2h^{\sigma}}z_{0}+\frac{d_{-,k}}{h^{\sigma}}\l[\l(1-\frac{\sigma}{2}\r)z_{1}+\frac{\sigma}{2}z_{2}\r]+\frac{c_{k}}{h},
& k=j-1, \hbox{ } \\
           \frac{d_{-,k}\sigma}{2h^{\sigma}}z_{0}+\frac{d_{+,k}}{h^{\sigma}}\l[\l(1-\frac{\sigma}{2}\r)z_{1}+\frac{\sigma}{2}z_{2}\r],
& k=j+1, \hbox{ } \\
           \frac{d_{+,k}}{h^{\sigma}}\l[\l(1-\frac{\sigma}{2}\r)z_{k-j}+\frac{\sigma}{2}z_{k-j+1}\r],
& k>j+1, \hbox{ } \\
           \frac{d_{-,k}}{h^{\sigma}}\l[\l(1-\frac{\sigma}{2}\r)z_{j-k}+\frac{\sigma}{2}z_{j-k+1}\r],
& k<j-1. \hbox{ }
         \end{array}
       \right.  
\eeee
Consequently, the formal solution of (\ref{e2-7}) can be written as 
$$v(t+\tau)=e^{\tau S}v(t)+\int_{t}^{t+\tau}e^{(t+\tau-\xi)S}p(v)d\xi,~~t,\tau\geq0.$$

Let $\tau$ be sufficiently small. While a trapezoidal rule is adopted for evaluating the integral,
we approximate the matrix exponential functions by the [1/1] Pad\'{e} approximant, that is,
$$e^{\tau S}=\l(I-\frac{\tau}{2}S\r)^{-1}\l(I+\frac{\tau}{2}S\r)+\O(\tau^{2}), ~~\tau\rightarrow0^{+},$$
in the formal solution.

If the step $\tau$ is determined via an adaptive procedure at each temporal level, such as the arc-length formula used in
\cite{zhu2020note,sk2001,Josh3},
then we acquire following fully discretized semi-adaptive scheme from (\ref{e2-7}) based on the initial value $v^0:$
\bb{e2-9}
v^{n}=\l(I-\frac{\tau_{n}}{2}S\r)^{-1}\l(I+\frac{\tau_{n}}{2}S\r)\l(v^{n-1}+\frac{\tau_{n}}{2}p^{n-1}\r)+\frac{\tau_{n}}{2}p^{n},
~~n=1,2,\cdots,M.
\ee

\section{Theoretical analysis of numerical solutions}
\label{sec: sec3}
In this section, we prove the monotonicity, positivity and stability of numerical solutions obtained from (\ref{e2-9}).
To this end, we must explore certain characteristic features of the matrices $I-\frac{\tau_{n}}{2}S, I+\frac{\tau_{n}}{2}S$ first
in following lemmas. Denote
$(d_{+})_{\max} = \max_{0\leq k\leq L}d_{+}(x_k,t),~(d_{-})_{\max} = \max_{0\leq k\leq L}d_{-}(x_k,t),$
and $c_{\max}= \max_{0\leq k\leq L}c(x_k,t).$

\red
\begin{definition}\label{def1}\cite{torres2024,Henrici1962}
If a matrix $A=(a_{ij})$ satisfies the following properties: 1) $a_{i,j}\leq 0, i\neq j; a_{ii} >0;$ 2) strictly diagonally dominant;
3) the sum of each row element in the matrix is greater than zero, then the weak-row sum criterion follows and the matrix
can be called an $M$-matrix in numerical algebra. Then all elements of the matrix $A^{-1}$ are positive.
In such a case, The matrix $A$ is said to be inverse positive.
\end{definition}
\black

\begin{lemma}\label{lem1}
If $\frac{\sqrt{17}-1}{2}\leq \sigma\leq 2$ and $\frac{\tau_{n}}{h^{\sigma}}<\frac{1}{2[(d_{+})_{\max}+(d_{-})_{\max}]+c_{\max}},$
then the matrix $I-\frac{\tau_{n}}{2}S$ is strictly diagonally dominant for $0<\tau_n\ll 1.$
\end{lemma}
\proof
Denote $P=I-\frac{\tau_{n}}{2}S=\l[p_{k,j}\r]\in \R^{L\times L}.$ It follows that
\bbbb
p_{k,k}&=&1+\frac{\tau_{n}(d_{+,k}+d_{-,k})}{2h^{\sigma}}\l(\frac{\sigma^{2}}{2}+\frac{\sigma}{2}-1\r)+\frac{\tau_{n}c_{k}}{2h};\\
p_{k,k+1}&=&-\frac{\tau_{n}d_{-,k}\sigma}{2h^{\sigma}}\l(\frac{\sigma^{2}+\sigma-4}{4}\r)-\frac{\tau_{n}d_{+,k}\sigma}{4h^{\sigma}}-\frac{\tau_{n}c_{k}}{2h};\\
p_{k,k-1}&=&-\frac{\tau_{n}d_{-,k}\sigma}{4h^{\sigma}}-\frac{\tau_{n}d_{+,k}\sigma}{2h^{\sigma}}\l(\frac{\sigma^{2}+\sigma-4}{4}\r);\\
p_{k,j}&=&-\frac{\tau_{n}d_{+,k}}{2h^{\sigma}}\l[\l(1-\frac{\sigma}{2}\r)z_{k-j}+\frac{\sigma}{2}z_{k-j+1}\r],~1\leq j < k-1;\\
p_{k,j}&=&-\frac{\tau_{n}d_{-,k}}{2h^{\sigma}}\l[\l(1-\frac{\sigma}{2}\r)z_{j-k}+\frac{\sigma}{2}z_{j-k+1}\r],~k+1< j \leq L.
\eeee
If $\frac{\sqrt{17}-1}{2}\leq\sigma\leq2,$ then
\bbbb
\frac{\sigma}{2}z_{k-j+1}+\l(1-\frac{\sigma}{2}\r)z_{k-j}>0,&&k>j,~k=2,3,\cdots,L;\\
\frac{1}{4}\l(\sigma^{2}+\sigma-4\r)\geq0,&&\frac{\sigma^{2}}{2}+\frac{\sigma}{2}-1>0.
\eeee
Therefore we have
$p_{k,k}>0,~p_{k,k-1}<0,~p_{k,k+1}<0,$ and $p_{k,j}<0,$ when $1\leq j<k-1$ or $k+1<j\leq L.$
Furthermore, sums of the absolute values of each row elements, except the diagonal element, of $P$
can be calculated via

\bbbb
r_{k}&=&\sum_{j=1,j\neq k}^{L}|p_{k,j}|
~=~\frac{\tau_{n}}{2}\cdot\frac{d_{+,k}}{h^{\sigma}}
\l\{\sum_{j=1}^{k-1}\l[\l(1-\frac{\sigma}{2}\r)z_{k-j}+\frac{\sigma}{2}z_{k-j+1}\r]+\frac{\sigma}{2}z_{0}\r\}\\
&&+\frac{\tau_{n}}{2} \cdot \frac{d_{-,k}}{h^{\sigma}}\l\{\sum_{j=k+1}^{L}\l[\l(1-\frac{\sigma}{2}\r)z_{j-k}
+\frac{\sigma}{2}z_{j-k+1}\r]+\frac{\sigma}{2}z_{0}\r\}+\frac{\tau_{n}}{2}\cdot\frac{c_{k}}{h}, ~1\leq k< L.
\eeee

Now, recall that 
$$\sum_{m=0}^{M}z_{m}\leq 0,~~M\in \Z^{+}.$$
Further, 
\bbbb
r_{k}&=&\sum_{j=1,j\neq k}^{L}|p_{k,j}|~\leq~\frac{\tau_{n}c_{k}}{2h}+\frac{\tau_{n}(d_{+,k}+d_{-,k})}
{2h^{\sigma}}\l(\frac{\sigma^{2}}{2}+\frac{\sigma}{2}-1\r)\\
&<&\frac{\tau_{n}c_{k}}{2h}+\frac{\tau_{n}(d_{+,k}+d_{-,k})}{2h^{\sigma}}
\l(\frac{\sigma^{2}}{2}+\frac{\sigma}{2}-1\r)+1~=~p_{k,k},~~k=1,2,\ldots,L.
\eeee
Moreover, we have
$p_{k,k}>0,~p_{k,k-1}<0,~p_{k,k+1}<0,$ and $p_{k,j}<0$ for $1\leq j<k-1$ or $k+1<j\leq L.$
Therefore $P$ is strictly diagonally dominant and, consequently, it is nonsingular.
\qedsymbol

\begin{lemma}\label{lem2}
If $\frac{\sqrt{17}-1}{2}\leq \sigma\leq 2$ and $\frac{\tau_{n}}{h^{\sigma}}<\frac{1}{2[(d_{+})_{\max}+(d_{-})_{\max}]+c_{\max}}$
then the matrix $I-\frac{\tau_{n}}{2}S$ is strictly monotone and inverse positive for $0<\tau_n\ll 1.$
\end{lemma}
\proof
Continuing from the proof of Lemma 3.1, we observe that
\bbbb
r_{k}&=&\sum_{j=1,j\neq k}^{L}|p_{k,j}|\leq\frac{\tau_{n}c_{k}}{2h}+\frac{\tau_{n}(d_{+,k}+d_{-,k})}
{2h^{\sigma}}\l(\frac{\sigma^{2}}{2}+\frac{\sigma}{2}-1\r)\\
&<&\frac{\tau_{n}c_{k}}{2h}+\frac{\tau_{n}(d_{+,k}+d_{-,k})}{2h^{\sigma}}
\l(\frac{\sigma^{2}}{2}+\frac{\sigma}{2}-1\r)+1~=~p_{k,k},~~k=1,2,\ldots,L.
\eeee
Furthermore, for $k\neq L,$
\bbbb
\sum_{j=1}^{L}p_{k,j}
&=&1-\frac{\tau_{n}d_{+,k}}{2h^{\sigma}}\l[\l(1-\frac{\sigma}{2}\r)\sum_{j=0}^{k-1}z_{j}
+\frac{\sigma}{2}\sum_{j=0}^{k}z_{j}\r]\\
&&-\frac{\tau_{n}d_{-,k}}{2h^{\sigma}}\l[\l(1-\frac{\sigma}{2}\r)\sum_{j=0}^{L-k}z_{j}
+\frac{\sigma}{2}\sum_{j=0}^{L-k+1}z_{j}\r]~>~0.
\eeee
If $k=L$ then
\bbbb
\sum_{j=1}^{L}p_{L,j}
&=&1-\frac{\tau_{n}d_{+,L}}{2h^{\sigma}}\l[\l(1-\frac{\sigma}{2}\r)\sum_{j=0}^{L-1}z_{j}
+\frac{\sigma}{2}\sum_{j=1}^{L}z_{j}\r]\\
&&-\frac{\tau_{n}d_{-,L}}{2h^{\sigma}}\l[\l(1-\frac{\sigma}{2}\r)z_{0}
+\frac{\sigma}{2}z_{1}\r]-\frac{\tau_{n}d_{-,L}\sigma}{4h^{\sigma}}z_{0}+\frac{\tau_{n}c_{k}}{2h}.
\eeee
Based on (\ref{e2-4}) and inequalities given in the theorem, we must have $\sum_{j=1}^{L}p_{L,j}>0.$
Consequently, $P$ is strictly monotone and inverse positive \cite{Josh3,Henrici1962}. 
\qedsymbol

\begin{lemma}\label{lem3}
If $\frac{\sqrt{17}-1}{2}\leq \sigma\leq 2$ and $\frac{\tau_{n}}{h^{\sigma}}<\frac{1}{2[(d_{+})_{\max}+(d_{-})_{\max}]+c_{\max}},$
then the matrix $I+\frac{\tau_{n}}{2}S$ is positive and nonsingular for $0<\tau_n\ll 1.$
\end{lemma}
\proof
Let $Q=I+\frac{\tau_{n}}{2}S=[q_{i,j}]\in \R^{L\times L}.$ Hence,
$$q_{k,k} = 1-\frac{\tau_{n}(d_{+,k}+d_{-,k})}{2h^{\sigma}}\l(\frac{\sigma^{2}}{2}
+\frac{\sigma}{2}-1\r)-\frac{\tau_{n}c_{k}}{2h},~~k=1,2,\ldots,L.$$
Similar to the analysis of the signs of the elements in matrix $P$ in Lemma \ref{lem1}, we may see that all elements
of the matrix $Q$ are greater than zero, hence the matrix $Q$ is positive
\cite{sheng2023nonconventional,sheng1999compound,Henrici1962}. Moreover,
\bbbb
&&\l\|\frac{\tau_{n}}{2}S\r\|_{\infty}~=~\frac{\tau_{n}}{2}\max_{1\leq k\leq L}\l\{\sum_{j=1}^{L}|s_{k,j}|\r\}\\
&&~~~~~=~\frac{\tau_{n}}{2}\max_{1\leq k\leq L}\frac{d_{+,k}}{h^{\sigma}}
\l\{\sum_{j=1}^{k-1}\l[\l(1-\frac{\sigma}{2}\r)z_{k-j}+\frac{\sigma}{2}z_{k-j+1}\r]+\frac{\sigma}{2}
z_{0}-\l[\l(1-\frac{\sigma}{2}\r)z_{0}+\frac{\sigma}{2}z_{1}\r]\r\}\\
&&~~~~~~~+\frac{\tau_{n}}{2}\max_{1\leq i\leq L}\frac{d_{-,k}}{h^{\sigma}}
\l\{\sum_{j=k+1}^{L}\l[\l(1-\frac{\sigma}{2}\r)z_{j-k}+\frac{\sigma}{2}z_{j-k+1}\r]+\frac{\sigma}{2}
z_{0}-\l[\l(1-\frac{\sigma}{2}\r)z_{0}+\frac{\sigma}{2}z_{1}\r]\r\}\\
&&~~~~~~~+\frac{\tau_{n}}{2}\max_{1\leq k\leq L}\frac{2c_{k}}{h}
~\leq~\frac{\tau_{n}}{h^{\sigma}}[(d_{+})_{\max}+(d_{-})_{\max}]\l(\frac{\sigma^{2}}{2}+\frac{\sigma}{2}-1\r)
+\frac{\tau_{n}}{h^{\sigma}}c_{\max}\\
&&~~~~~\leq~\frac{\tau_{n}}{h^{\sigma}}\{2[(d_{+})_{\max}+(d_{-})_{\max}]+c_{\max}\}.
\eeee
Thus, to ensure  
$$\l\|\frac{\tau_{n}}{2}S\r\|_{\infty}<1,$$
a sufficient condition has to be
$$\frac{\tau_{n}}{h^{\sigma}}<\frac{1}{2[(d_{+})_{\max}+(d_{-})_{\max}]+c_{\max}}.$$
Therefore we must need constraints
$\frac{\tau_{n}}{h^{\sigma}}<\frac{1}{2[(d_{+})_{\max}+(d_{-})_{\max}]+c_{\max}}$
together with $\frac{\sqrt{17}-1}{2}\leq \sigma\leq 2.$ 
These ensure the nonsingularity of $Q$ and the lemma is clear.
\qedsymbol

~ 

Based on Lemmas \ref{lem1}-\ref{lem3}, we proceed for the stability proof of the semi-adaptive
finite difference scheme (\ref{e2-9}). We need:

\begin{definition}\label{def2}
Let $v^n,\widetilde{v}^n$ denote the exact and perturbed solutions of a finite
difference method such as {\rm(\ref{e2-9})}, respectively. Set $E^n=v^n-\widetilde{v}^n,~0\leq n\leq M.$
Then the numerical method is conditionally stable if there exists a positive constant $c$
such that $\|E^n\|\leq c\|E^0\|,~0\leq n\leq M,$ under certain restraints.
\end{definition}

\begin{theorem}\label{thm1}
Suppose $\frac{\sqrt{17}-1}{2}\leq \sigma\leq 2$ and $\frac{\tau_{n}}{h^{\sigma}}<\frac{1}{2[(d_{+})_{\max}+(d_{-})_{\max}]+c_{\max}}$
hold for all $n\geq 0.$ If the nonlinear term $p$ is frozen, that is,
$p^k\equiv \tilde{p}^k,~k=n,n+1,$
then the fully discretized semi-adaptive scheme $(\ref{e2-9})$ is conditionally stable.
\end{theorem}
\proof
Recall Lemmas discussed and $(\ref{e2-9}).$
Denote $\xi_{k}^{n}=\frac{\tau_{n}d_{+,k}}{2h^{\sigma}},\,\omega_{k}^{n}=\frac{\tau_{n}d_{-,k}}{2h^{\sigma}},\,
\varsigma_{k}^{n}=\frac{\tau_{n}c_{k}}{2h},$ $k=1,2,\ldots,L, n=0,1,\ldots,M.$
If the source term in scheme (\ref{e2-9}) is frozen then the following equalities become true: \black
$$E^{n+1}=\l(I-\frac{\tau_n}{2}S\r)^{-1}\l(I+\frac{\tau_{n}}{2}S\r){E^{n}},~~~n=0,1,\ldots,M-1.$$

We now proceed to prove our theorem through an induction. Assume that $\left|E_{l}^{1}\right|$ be
the maximal value of $|E_{k}^{1}|,~1\leq k\leq L.$ Firstly, when $n = 0,$ we have
\bbbb
\l\|E^{1}\r\|_{\infty}
&=&\max\limits_{1\leq k\leq L}|E_{k}^{1}| = |E_{l}^{1}|
\leq |E_{l}^{1}|-\xi_{l}^{1}\l(1-\frac{\sigma}{2}\r)\sum_{m=0}^{l-1}z_{m}|E_{l}^{1}|-\xi_{l}^{1}\frac{\sigma}{2}\sum_{m=0}^{l}z_{m}|E_{l}^{1}|\\
&&-~\omega_{l}^{1}\l(1-\frac{\sigma}{2}\r)\sum_{m=0}^{L-l}z_{m}|E_{l}^{1}|-\omega_{l}^{1}\frac{\sigma}{2}\sum_{m=0}^{L-l+1}z_{m}|E_{l}^{1}|\\
&\leq&\l[1+\xi_{l}^{0}\l(1-\frac{\sigma}{2}\r)\sum_{m=0}^{l-1}z_{m}+\frac{\sigma}{2}\xi_{l}^{0}\sum_{m=0}^{l}z_{m}
+\omega_{l}^{0}(1-\frac{\sigma}{2})\sum_{m=0}^{L-l}z_{m}\r.\\
&&+\l.\omega_{l}^{0}\frac{\sigma}{2}\sum_{m=0}^{L-l+1}z_{m}\r]\cdot\max\limits_{1\leq k\leq L}|E_{k}^{0}|\\
&\leq& \max\limits_{1\leq k\leq L}\l|E_{k}^{0}\r|=\l\|E^{0}\r\|_{\infty}. 
\eeee

When $\|E^{1}\|_{\infty}=|E_{L}^{1}|,$ it follows from properties of $z_{m}$ that
\bbbb
\|E^{1}\|_{\infty} &=& |E_{L}^{1}|
~\leq~ |E_{L}^{1}|-\xi_{L}^{1}\l(1-\frac{\sigma}{2}\r)\sum_{m=0}^{L-1}z_{m}|E_{L}^{1}|-\xi_{L}^{1}\frac{\sigma}{2}\sum_{m=1}^{L}z_{m}|E_{L}^{1}|\nnn\\
&&-~\omega_L^{1}\frac{\sigma}{2}z_{0}\l|E_{L-1}^{1}\r|-\omega_L^{1}\l[\l(1-\frac{\sigma}{2}\r)z_0+\frac{\sigma}{2}z_1\r]|E_L^{1}|\nnn\\
&\leq&\l[1+\zeta_{L}^{1}-(\xi_{L}^{1}+\omega_{L}^{1})\l(1-\frac{\sigma}{2}\r)z_{0}-(\xi_{L}^{1}+\omega_{L}^{1})\frac{\sigma}{2}z_{1}\r]|E_{L}^{1}|\nnn\\
&&-~\xi_{L}^{1}\sum_{m=2}^{L}\l[\l(1-\frac{\sigma}{2}\r)z_{m-1}+\frac{\sigma}{2}z_{m}\r]|E_{L-m+1}^{1}|-\omega_{L}^{1}\frac{\sigma}{2}z_{0}|E_{L-1}^{1}|\nnn\\
&\leq&\max\limits_{1\leq k\leq L}\l|E_{k}^{0}\r|=\l\|E^{0}\r\|_{\infty}. 
\eeee

Secondly, assume that $\l\|E^{\ell}\r\|_{\infty} \leq \l\|E^{0}\r\|_{\infty}$ for $\ell=1,2,\ldots, n-1,$ and denote $|E_{l}^{n}|$ as the maximum
of $|E_{k}^{n}|,~~1\leq k\leq L.$ We have
\bbbb
\|E^{n}\|_{\infty}
&=&\max\limits_{1\leq k\leq L}|E_{k}^{n}| = |E_{l}^{n}|
\leq|E_{l}^{n}|-\xi_{l}^{n}\l(1-\frac{\sigma}{2}\r)\sum_{m=0}^{l-1}z_{m}|E_{l}^{n}|
-\xi_{l}^{n}\frac{\sigma}{2}\sum_{m=0}^{l}z_{m}|E_{l}^{n}|\nnn\\
&&-~\omega_{l}^{n}\l(1-\frac{\sigma}{2}\r)\sum_{m=0}^{L-l}z_{m}|E_{l}^{n}|-\omega_{l}^{n}\frac{\sigma}{2}\sum_{m=0}^{L-l+1}z_{m}|E_{l}^{n}|\nnn\\
&\leq&\l[1+\zeta_{l}^{n}-(\xi_{l}^{n}+\omega_{l}^{n})\l(1-\frac{\sigma}{2}\r)z_{0}-(\xi_{l}^{n}+\omega_{l}^{n})\frac{\sigma}{2}z_{1}\r]|E_{l}^{n}|
-\xi_{l}^{n}\frac{\sigma}{2}z_{0}|E_{l+1}^{n}|\nnn\\
&&-~\zeta_{l}^{n}|E_{l+1}^{n}| -\xi_{l}^{n}\sum_{m=2}^{l}\l[\l(1-\frac{\sigma}{2}\r)z_{m-1}+\frac{\sigma}{2}z_{m}\r]|E_{l-m+1}^{n}|-\omega_{l}^{n}\frac{\sigma}{2}z_{0}|E_{l-1}^{n}|\nnn\\
&&-~\omega_{l}^{n}\sum_{m=2}^{L-l+1}\l[\l(1-\frac{\sigma}{2}\r)z_{m-1}+\frac{\sigma}{2}z_{m}\r]|E_{l+m-1}^{n}|\nnn\\
&\leq&\max\limits_{1\leq k\leq L}|E_{k}^{n-1}|=\l\|E^{n-1}\r\|_{\infty}.
\eeee

Note that the above inequalities hold for all $n\in\{1,2,\cdots,M\}.$ Therefore
$\l\|E^{n}\r\|_{\infty}\leq\l\|E^{0}\r\|_{\infty}$ is warranted for all $n\in\{1,2,\cdots,M\}.$ Likewise,
$\l\|E^{n}\r\|_{\infty}\leq \l|E_L^0\r|$ is secured since $\l\|E^{n}\r\|_{\infty}=\l|E_{L}^{n}\r|$ are true for
$n=0,1,\ldots,M-1.$
\qedsymbol

~

Next, we wish to show the positivity and monotonicity of the numerical solution generated by (\ref{e2-9}).
To this end, we need:
\begin{lemma}\label{lem4}
Let $v^{0}=0$ and $\tau_1\geq 0$ be an initial time-step. If $\frac{\sqrt{17}-1}{2}\leq \sigma\leq 2$, and
$\frac{\tau_{1}}{h^{\sigma}}<\frac{1}{2[(d_{+})_{\max}+(d_{-})_{\max}]+c_{\max}}$ holds together with
$\tau_1<1/ \varrho,$ where $\varrho=\max(p(\tau_1 p_{0})),$ then $v^{1}>v^{0}$ and $\l\|v^1\r\|_{\infty}<1.$
\end{lemma}
\proof
Exploiting (\ref{e2-9}), we observe that
\bbb
v^{1}&=& \l(I-\frac{\tau_1}{2}S\r)^{-1}\l(I+\frac{\tau_1}{2}S\r)\l(v^{0}+\frac{\tau_1}{2}{p}^{0}\r) + \frac{\tau_1}{2}{p}^{1}\nnn\\
&=&\frac{\tau_1}{2}\l[\l(I-\frac{\tau_1}{2}S\r)^{-1}\l(I+\frac{\tau_1}{2}S\r){p}^{0}+{p}(\tau_1 {p}^{0})\r]. \label{e3-4}
\eee

It follows immediately that  $v^{1}>v^{0}=0$ as a consequence of Lemmas \ref{lem2} and \ref{lem3}.
Next, we may demonstrate that $\l\|v^{1}\r\|_{\infty}<1.$ For this, we define $Y=(1,1,\ldots,1)^{\tT}$. Then,
$$v^{1}-Y=\l(I-\frac{\tau_{1}}{2}S\r)^{-1}\l[\l(\frac{\tau_{1}}{2}S+I\r)\frac{\tau_{1}}{2}p^{0}+\l(I-\frac{\tau_{1}}{2}S\r)
\frac{\tau_{1}}{2}p(\tau_{1}p^{0})-\l(I-\frac{\tau_{1}}{2}S\r)Y\r].$$
We further define
\bbbb
Z_{1}&=&\l(\frac{\tau_{1}}{2}S+I\r)\frac{\tau_{1}}{2}p^{0}+\l(I-\frac{\tau_{1}}{2}S\r)\frac{\tau_{1}}{2}p(\tau_{1}p^{0}),~
Z_{2}~=~-\l(I-\frac{\tau_{1}}{2}S\r)Y,\\
Z&=&Z_{1}+Z_{2}.
\eeee
It may be observed readily that
\bbb
Z_{1}
&\leq&\l[\l(\frac{\tau_{1}}{2}S+I\r)+\l(I-\frac{\tau_{1}}{2}S\r)\r]\frac{\tau_{1}}{2}p(\tau_{1}p^{0})=\tau_{1}p(\tau_{1}p^{0})\label{e3-6},\\
Z_{2}&=& -Y+\frac{\tau_1}{2}SY \leq -Y.\nnn
\eee
Hence,
$$Z=Z_{1}+Z_{2}\leq \tau_{1}\varrho-1.$$
Apply the above results we acquire that
$$v^{1}-Y=\l(I-\frac{\tau_{1}}{2}S\r)\l(Z_{1}+Z_{2}\r)\leq\l(I-\frac{\tau_{1}}{2}S\r)\l(\tau_{1}\varrho-1\r).$$

Combining the above with Lemma \ref{lem1}, we have $v^{1}<Y$ if $\tau_1<1/\varrho.$ Therefore $\|v^{1}\|_{\infty}<1.$
This completes the proof.
\qedsymbol

\begin{theorem}\label{thm2}
If $\frac{\sqrt{17}-1}{2}\leq \sigma\leq 2$ and $\frac{\tau_{n}}{h^{\sigma}}<\frac{1}{2[(d_{+})_{\max}+(d_{-})_{\max}]+c_{\max}}$
holds for all $\tau_{n},$ and $0\leq v^{n-1} <1$ such that $Sv^{n-1}+{p}^{n-1}>0,~n\geq 0,$ then
\begin{enumerate}
\item[\rm(1)] the solution sequence generated by {\rm(\ref{e2-9})}, $\l\{v^{n}\r\}_{n=0}^{\infty},$ is strictly monotonically increasing;
\item[\rm(2)] $ Sv^{n}+{p}^{n}>0,~n\geq 0.$
\end{enumerate}
\end{theorem}
\proof
To begin, we recall Lemma \ref{lem4} and the inequality $v^{1}>v^{0}.$ Next, similar to strategies used in
\cite{zhu2020note,padgett2018quenching,martinez2022nonlinear,liu2023semi}, we approximate
$p^{n}$ in (\ref{e2-9}) through an Euler's formula, that is,
\bb{e3-7}
{p}^{n}=p^{n-1}+\tau_n M (Sv^{n-1}+p^{n-1})+\O(\tau_n^2),~~n\geq 1,
\ee
where $M$ is the nonnegative diagonal Jacobian matrix of $p(v).$
Recall (\ref{e2-9}). We find that
\bbbb
v^{n}-v^{n-1}&=& \l(I-\frac{\tau_{n}}{2}S\r)^{-1}\l[\l(I+\frac{\tau_{n}}{2}S\r)\l(v^{n-1}+\frac{\tau_{n}}{2}{p}^{n-1}\r) \r.\nnn\\
&&\l. + \frac{\tau_{n}}{2}\l(I-\frac{\tau_{n}}{2}S\r){p}^{n} - \l(I-\frac{\tau_{n}}{2}S\r)v^{n-1}\r]\nnn\\
&=& \l(I-\frac{\tau_{n}}{2}S\r)^{-1}\l[\tau_{n}S v^{n-1}+\frac{\tau_{n}}{2}{p}^{n-1}
+\frac{\tau_{n}}{2}{p}^{n}+\frac{\tau_{n}^{2}}{4}S({p}^{n-1}-{p}^{n})\r]\nnn\\
&\geq& \tau_{n}\l(I-\frac{\tau_{n}}{2}S\r)^{-1}\l[S v^{n-1}+{p}^{n-1}
- \frac{\tau_{n}}{4}S({p}^{n}-{p}^{n-1})\r],~n>1,
\eeee
due to the property of ${p}^{n}>{p}^{n-1}.$ Furthermore, 
$$v^{n}-v^{n-1}\geq \tau_{n}\l(I-\frac{\tau_{n}}{2}S\r)^{-1}\l[(Sv^{n-1}+{p}^{n-1})-\frac{M\tau_{n}^2}{4}S(Sv^{n-1}+{p}^{n-1})\r].$$

Hence, $v^{n}>v^{n-1}$ if $\tau_{n}$ is sufficiently small. Besides,
it can be seen that
\bbbb
S v^{n}+{p}^{n}&=&{p}^{n}-p^{n-1}+S v^{n-1}+{p}^{n-1}+S(v^{n}-v^{n-1})\nnn\\
&\geq& {p}^{n}+\l(I-\frac{\tau_{n}}{2}S\r)^{-1}\l[\l(I-\frac{\tau_{n}}{2}S\r)\l(S v^{n-1}+{p}^{n-1}\r)\r.\nnn\\
&&+\l.\tau_{n} S({p}^{n-1}+S v^{n-1}) + \frac{\tau_{n}^2}{4}S^{2}({p}^{n-1}-{p}^{n})\r]-{p}^{n-1}\nnn\\
&=&{p}^{n}-{p}^{n-1}+\l(I-\frac{\tau_{n}}{2}S\r)^{-1}\l[\l(I+\frac{\tau_{n}}{2}S\r)\l({p}^{n-1}+S v^{n-1}\r)+\frac{\tau_{n}^2}{4}S^{2}({p}^{n-1}-{p}^{n})\r],
\label{e3-8}
\eeee
since ${p}^{n}\approx {p}[v^{n-1}+\tau_{n} (Sv^{n-1}+{p}^{n-1})]$ and $p^{n}>p^{n-1}.$
 In accordance with Taylor's theorem, the following is true:
\bbbb
S v^{n}+{p}^{n}&\geq& 
{p}^{n}-{p}^{n-1}+\l(I+\tau_{n}S-\frac{\tau_{n}^{3}}{4}MS^2\r)(S v^{n-1}+{p}^{n-1}),
\eeee
where $M$ is the Jacobian matrix used in (\ref{e3-7}).   
Consequently, for any sufficiently small $\tau_{n},$ the matrix $I+\tau_{n}S-\frac{M \tau_{n}^{3}}{4}S^2$
must be positive \cite{Henrici1962}.  
Therefore our Theorem is proved.
\qedsymbol

~

In light of our investigation, we may state the following corollary for the vector solution
sequence generated by the semi-adaptive method (\ref{e2-9}).

\begin{corollary}\label{thm3}
If $\frac{\sqrt{17}-1}{2}\leq \sigma\leq 2$ and $\frac{\tau_{\ell}}{h^{\sigma}}<\frac{1}{2[(d_{+})_{\max}+(d_{-})_{\max}]+c_{\max}}$
hold for all $0\leq \ell\leq n,$ and $v^{0}\geq 0$ so that $Sv^{0}+{p}^{0}>0,$ then the solution sequence
$\left\{v^{\ell}\,\right\}_{\ell = 0}^{n}$ generated by the semi-adaptive finite difference scheme
$(\ref{e2-9})$ is strictly monotonically increasing until a quenching is reached.
\end{corollary}

\section{Simulations of fractional order quenching phenomena }
\label{sec: sec4}

Four simulation experiments in connection to the
one-dimensional, two-sided Riemann-Liouville fractional order convection-diffusion problem
(\ref{e2-1}) will be carried out and analyzed in this section. Comparisons with existing results
will be conducted. Impacts of the variable convective coefficient $c(x,t)={b}/{x}, 0<x<a,$
and fractional order $\sigma$ to
the critical interval size, quenching time, and quenching location will be
examined in details, respectively. A weak Courant constraint
\cite{sheng2023nonconventional,Josh3,mooney1996implicit} is incorporated
throughout computations.
\color{red} The numerical experiments are implemented as follows:
\begin{enumerate}
  \item Initialisation: Set up the spatial grid points $x_{k}=kh,~k=0,1,2,\ldots,L;$ the temporal grid points $t_{n}=n\tau,~n=0,1,2,\ldots,M.$ And the initial condition $v(x,0)=\psi(x)$.
  \item Discretization: Discretize the fractional order derivatives using standard and shifted weighted average of Gr\"{u}nwald formula and discretize the convective terms using Euler's formula,then (\ref{e2-1}) turn into a ordinary differential (2.6) . We construct a semi-discretized differential format (\ref{e2-9}).
  \item Update the solution: Solve the discretized system of nonlinear equations using an iterative method to obtain the numerical solution $v^{n}$.
\end{enumerate}
\color{black}

\subsection{Simulation experiment A: critical lengths for quenching}
\label{sec: sec4-1}

Consider (\ref{e2-1}), we adopt a basic pair of the temporal step $\tau=2\times10^{-4}$ and spacial step $h=a\times10^{-2}.$
Fix $\sigma=2.$ Further, set $\psi(x)=v(x,0)\equiv 0.$ The quenching critical interval size $a^{*}$
is then evaluated through the solution sequence from (\ref{e2-9}). They are then compared with
known results of $a^{*S}$ for integer-order problems summarized in \cite{sheng1999compound}.

\begin{table}
\begin{center}
 \tabcolsep 0.05in\small

\caption{Critical interval sizes $a^{*}$ vs. documented $a^{*S}$ values ($\sigma=2$). Different $b$ values are used.}
\vspace{2mm}

\begin{tabular}{lllllllll}
 \hline
  $b$ & $0.95$ & $0.9$ & $0.4$ & $0$ &$-0.4$   & $-0.9$ & $-1.2$ & $-2.0$\\\hline
  $a^{*}(\theta=2)$ & 1.0297 & 1.0312 & 1.0936 &1.1832 &1.2788  & 1.3914 & 1.4548 & 1.6123 \\
  $a^{*S}(\theta=2)$ & 1.1218 & 1.1218 & 1.1283 & $ $ & 1.2590  & 1.3576 & 1.4167 & 1.5657 \\
  $a^{*}(\theta=1)$ & 1.3304 & 1.3326 & 1.4149 & 1.5303 &1.6522  & 1.7949 & 1.8749 & 2.0726\\
  $a^{*S}(\theta=1)$ & 1.4479 & 1.4487 & 1.4788 & $ $ & 1.6097  & 1.7315 &1.8280 & 2.0124\\
  $a^{*}(\theta=0.5)$ & 1.6368 & 1.6401 & 1.7445 & 1.8856 & 2.0329 & 2.2038 & 2.2989  & 2.5327 \\
  $a^{*S}(\theta=0.5)$ & 1.7800 & 1.7800 & 1.8004 & $ $ & 2.0037  & 2.1521 & 2.2389 & 2.5000 \\[2pt]\hline
 \end{tabular}
 \label{tab1}
 \end{center}
\end{table}

\begin{figure}
\begin{center}
\epsfig{file=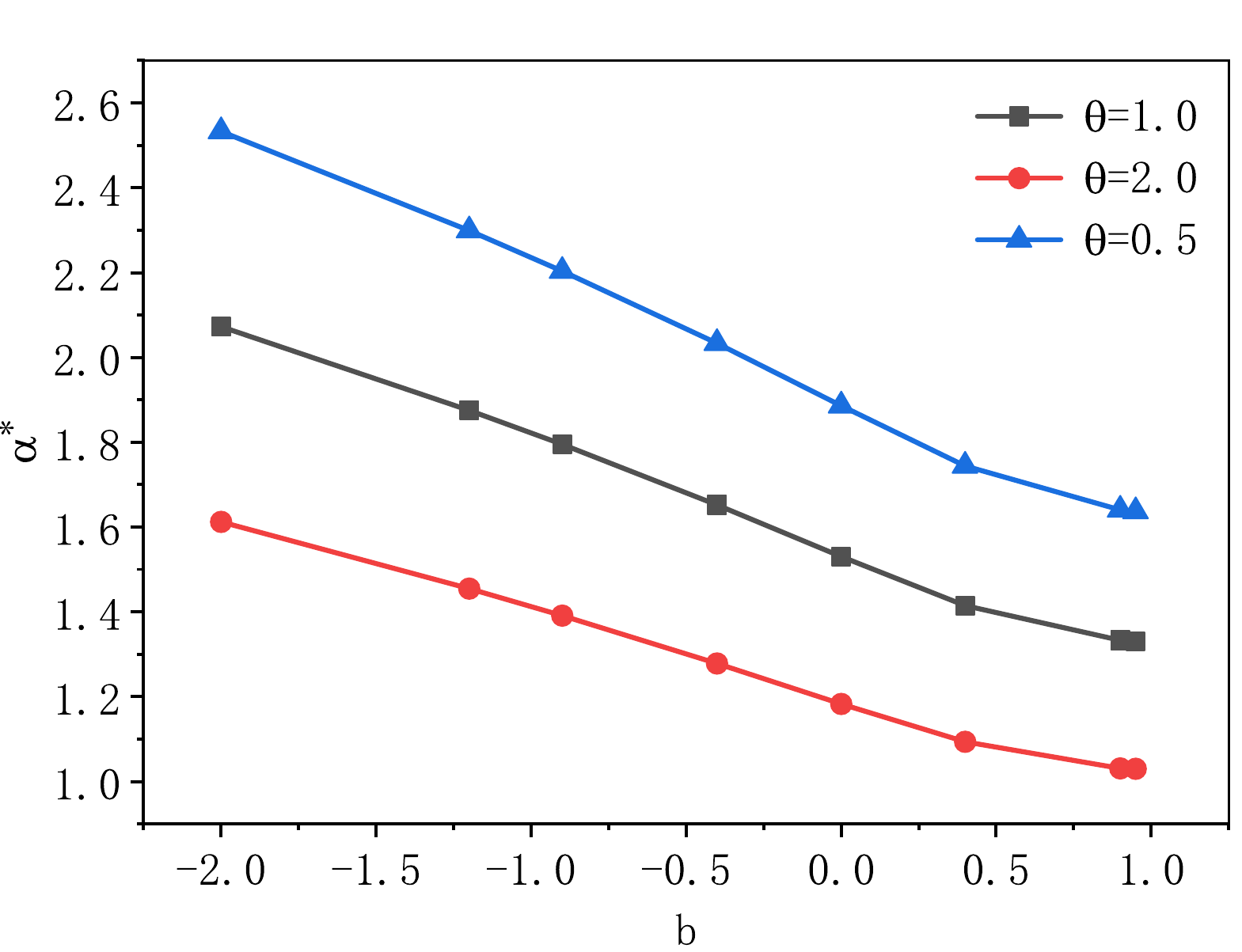,width=7.3cm,height=5.4cm}
\epsfig{file=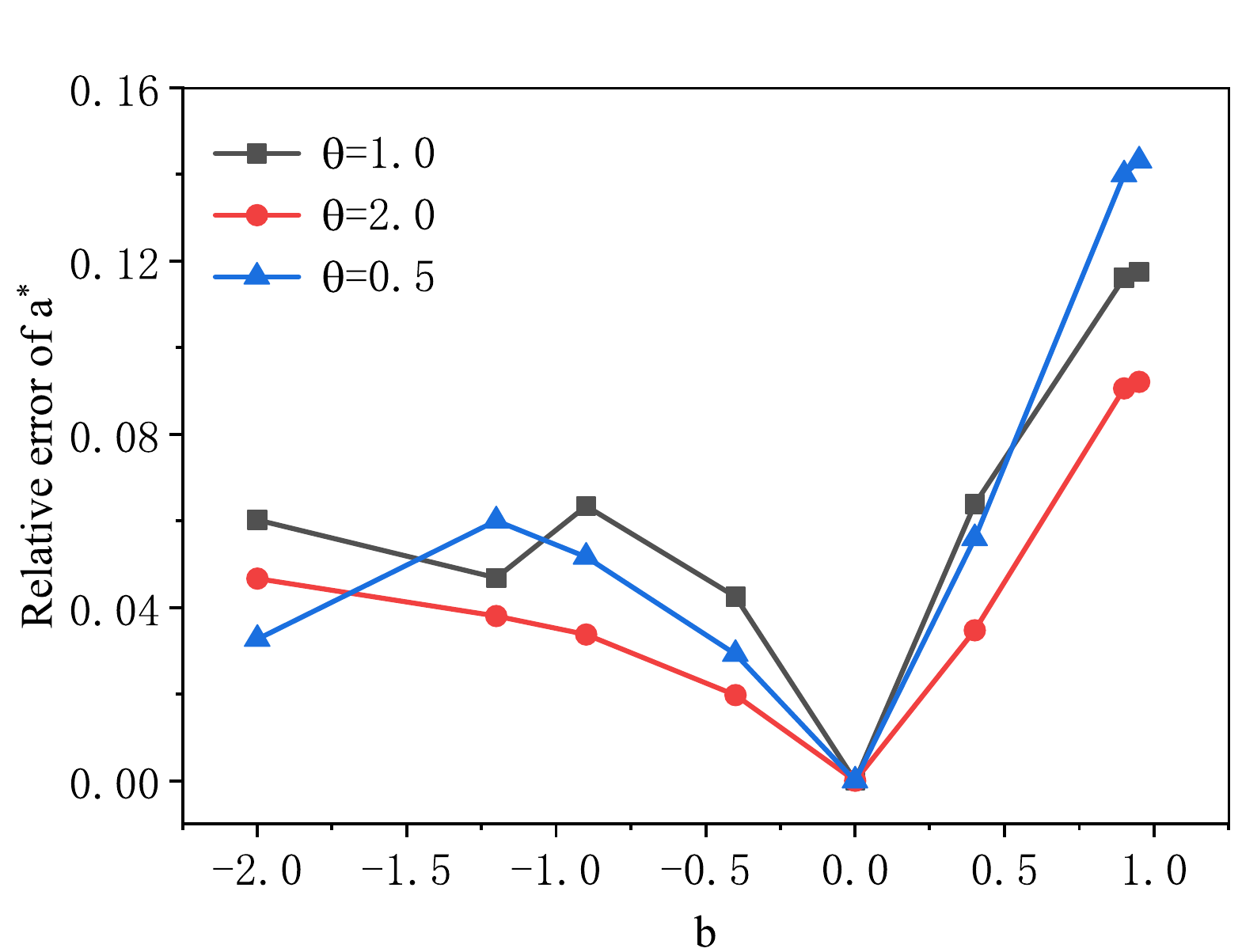,width=7.3cm,height=5.6cm}
\caption{LEFT:  Relationship between $a^{*}$ and  $b$ ($\sigma=2$);
RIGHT: Relationship between relative error $\left|a^{*}-a^{*S}\right|/\left|a^{*S}\right|$ and parameter $b.$ }  
\label{fig1}
\end{center}
\end{figure}

Numerical results of the critical size $a^{*}$ versus the convective constant $b$
corresponding to different $\theta$ values are listed in Table \ref{tab1}.
It is evident that $a^{*}$ increases as $b$ decreases.
The same phenomenon is illustrated in Fig. \ref{fig1} (left)
which simulates the relationship between the length $a^{*}$ and $b$ for $\theta=0.5, 1, 2,$ respectively.
On the other hand, Fig. \ref{fig1} (right) exhibits the relative error profile between $a^{*}$ found in this
paper and the standard estimates $a^{*S}$ in \cite{sheng1999compound}. Our experiments indicate that
$a^{*}$ decreases as $\theta$ increases, which is truly in line with known observations shown in \cite{sheng1999compound}.

\begin{table}
\begin{center}
 \tabcolsep 0.05in\small
\caption{Correlations between critical regions $a^{*}$ and order $\sigma$ of fractional derivative, coefficient $b.$}
\vspace{2mm}

\begin{tabular}{lllllllll}
 \hline
  $\sigma$ & $1.55$ & $1.6$ & $1.7$ & $1.8$ & $1.9$   & $2.0$\\\hline
  $a^{*},b=-2.0$ & 1.7500 & 1.7810 & 1.8612 & 1.9414 & 2.0126 & 2.0726 \\
  $a^{*},b=-1.2$ & 1.4759 & 1.5281 & 1.6292 & 1.7227 & 1.8053 & 1.8749 \\
  $a^{*},b=-0.9$ & 1.3667 & 1.4231 & 1.5320 & 1.6622 & 1.7205 & 1.7949 \\
  $a^{*},b=-0.4$ & 1.1571 & 1.2232 & 1.3501 & 1.4649 & 1.5660 & 1.6522 \\
  $a^{*},b=0$ & 0.9230 & 1.0072 & 1.1643  & 1.3047  & 1.4270 & 1.5303 \\
  $a^{*},b=0.4$ & 0.8492 & 0.9016 & 1.0237 & 1.1614 & 1.2963 & 1.4149 \\
  $a^{*},b=0.9$ & 1.0410 & 1.0638 & 1.1145 & 1.1728 & 1.2452 & 1.3326 \\
  $a^{*},b=0.95$ & 1.0609 & 1.0823 & 1.1296 & 1.1830 & 1.2488 & 1.3304 \\ [2pt]\hline
  \end{tabular}
  \label{tab2}
 \end{center}
\end{table}

%

\begin{figure}
\begin{center}
\epsfig{file=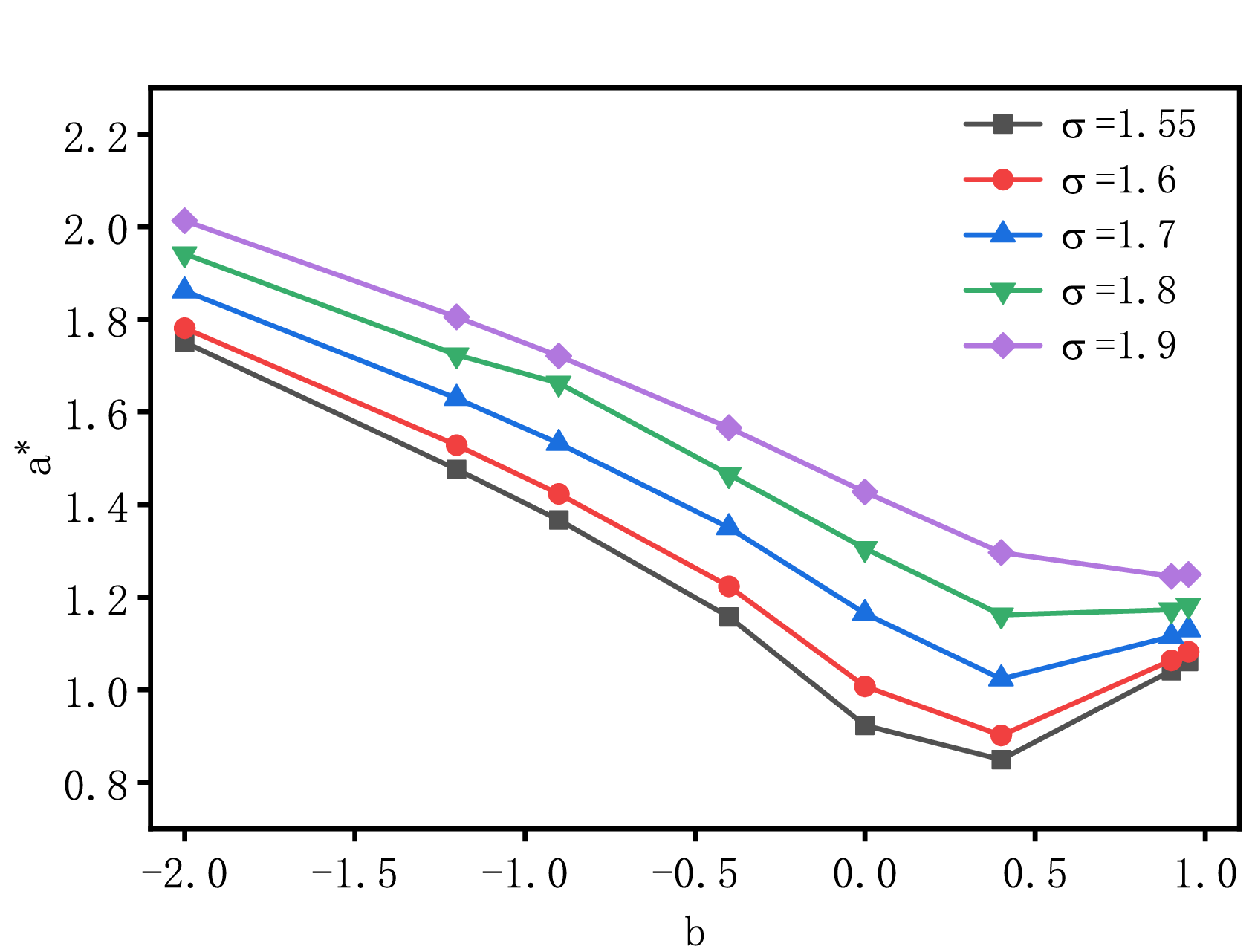,width=7.3cm,height=5.8cm}
\epsfig{file=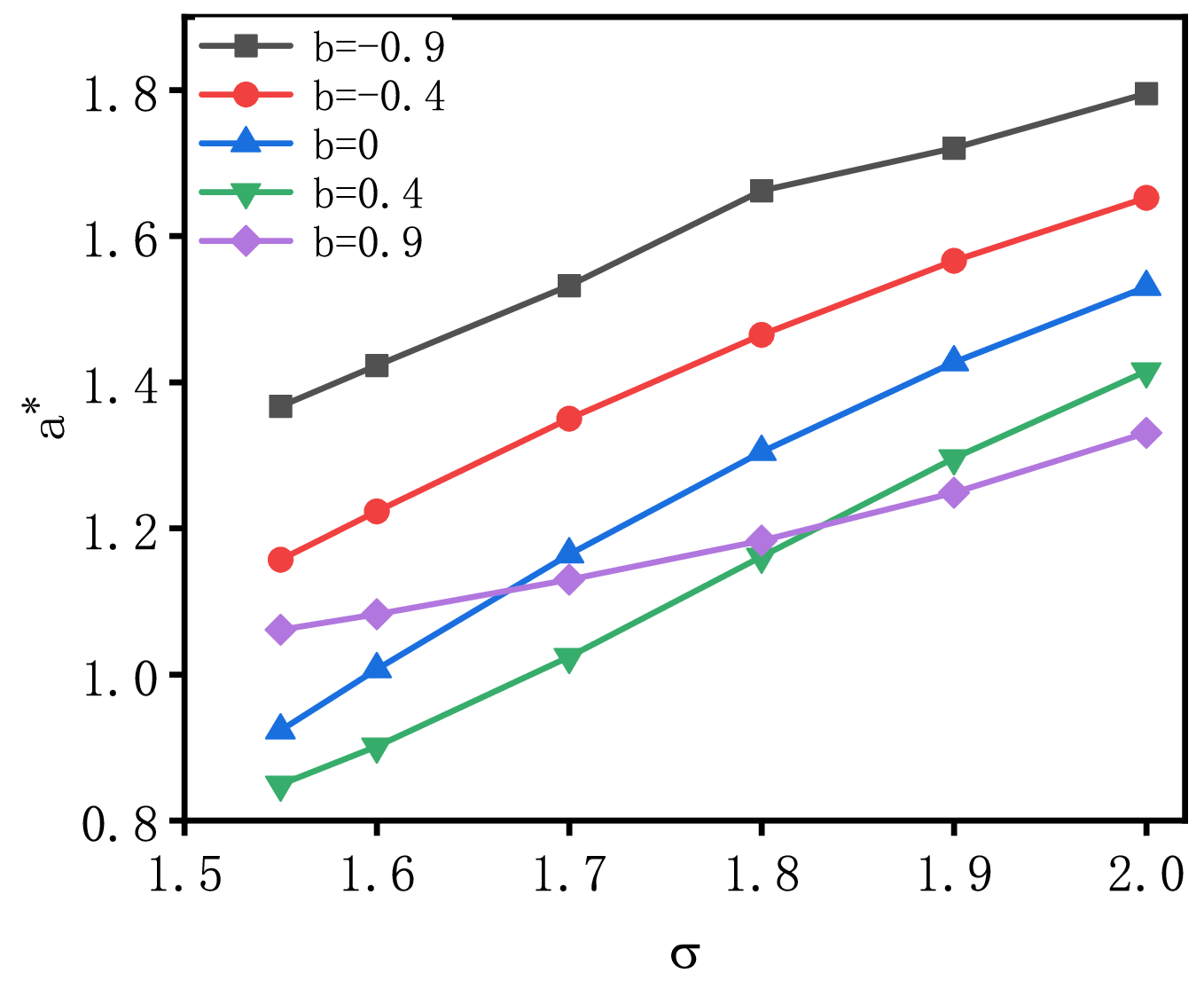,width=7.3cm,height=5.4cm}
\caption{LEFT: Relationship between $a^{*}$ and $b$ ($\theta=1$);
RIGHT: Relationship between $a^{*}$ and $\sigma$ ($\theta=1$).}
\label{fig2}
\end{center}
\end{figure}

Table \ref{tab2} illustrates the critical interval length $a^{*}$ in connection with the order of fractional derivatives
$\sigma$ based on different values of $b.$ We may observe that, though computed values of $a^{*}$ slowly diminish
as $b$ increases initially, they bounce back to increase from $b=0.4$ and beyond. Further, for a fixed value of $b,$
the critical interval length $a^{*}$ monotonically increases as the order $\sigma$ increases. Such monotone
trajectories of $a^{*}$ are more evident when $b=0$ as we reported earlier in \cite{zhu2023simulation}.


\subsection{Simulation experiment B: quenching time of the solutions}
\label{sec: sec4-2}

A motivation of this study is to explore further interconnections between the quenching time $T_a$ and
interval size $a,$ coefficient $b,$ and fractional order $\sigma.$ It becomes evident that $T$ can be affected
by $a,b,\sigma.$ This completely changes the traditional understanding that $T_a$ depends only on the
interval size $a$ \cite{sk2001,sheng1999compound}. The implicit scheme (\ref{e2-9}) is used throughout the exploration.

\begin{figure}
\begin{center}
\epsfig{file=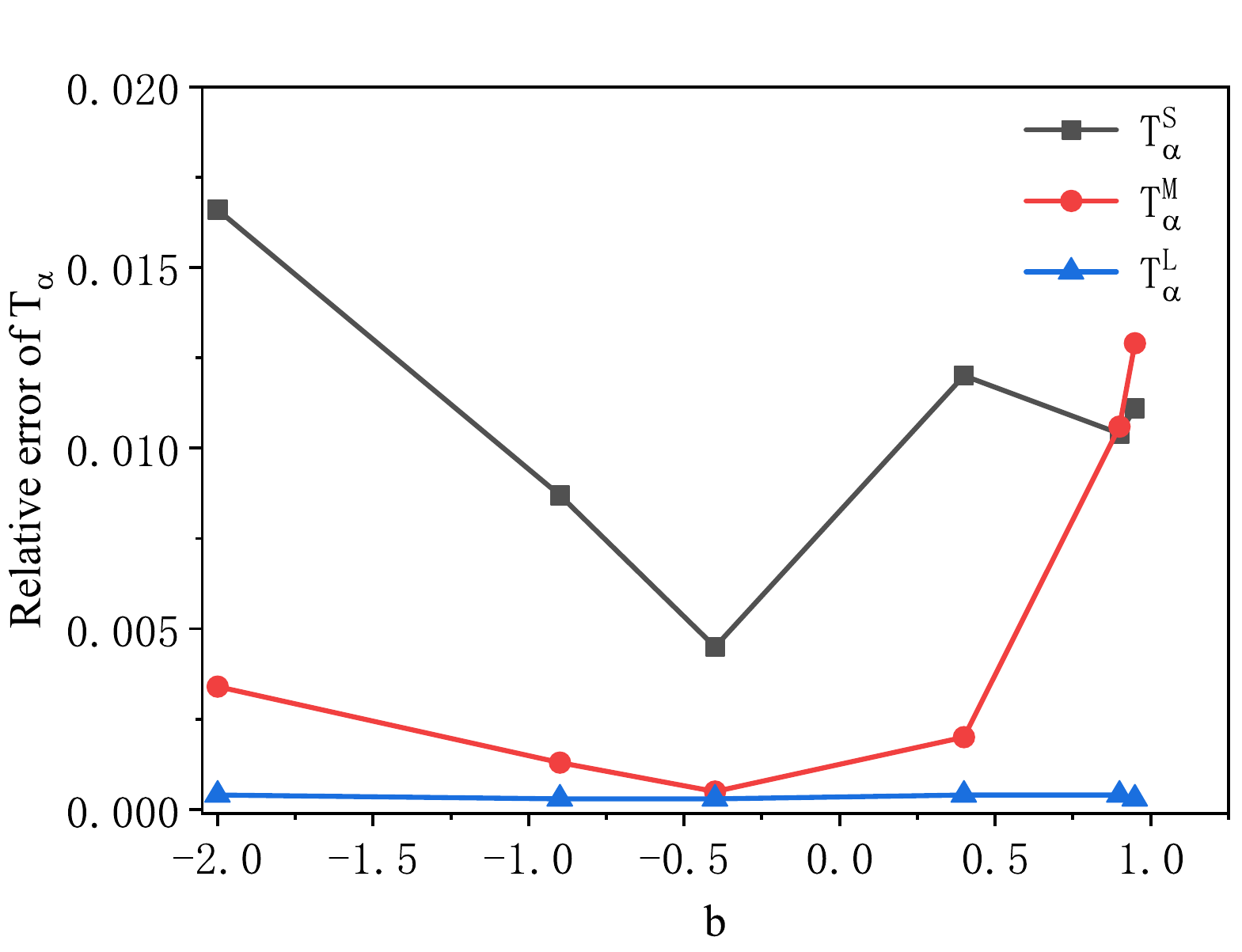,width=7.2cm,height=5.2cm}
\epsfig{file=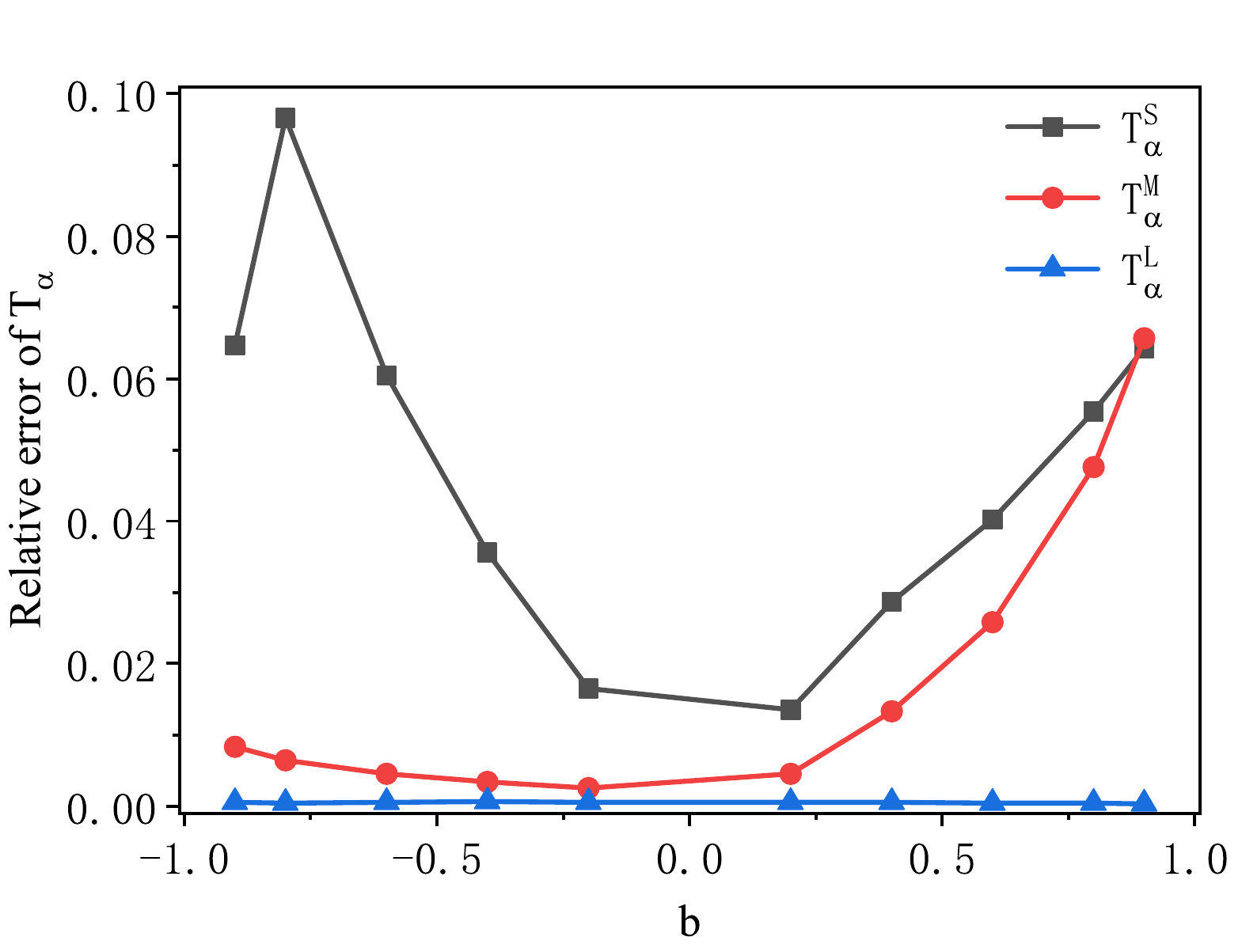,width=7.2cm,height=5.2cm}
\caption{LEFT: Relative error of computed quenching time $T_{a}$ from $T_{a}^{S},T_{a}^{M},T_{a}^{L}$
($a=\pi,\sigma=2,\theta=1$);  RIGHT: Relative error of computed quenching time $T_{a}$ from $T_{a}^{S},T_{a}^{M},T_{a}^{L}$
($a=2,\sigma=2,\theta=1$).}
\label{fig3}
\end{center}
\end{figure}

Firstly, in Tables \ref{tab3} and \ref{tab4}, we show the quenching time $T_{a}$ as $b$ varies.
Secondly, we compare our results obtained with experimental data reported by Sheng
\cite{sheng1999compound}, Liu \cite{liu2023semi} and Mooney \cite{mooney1996implicit} et al. (values
are denoted $T_{a}^{S},\,T_{a}^{L},\,T_{a}^{M},$ respectively).
Fig. \ref{fig3} presents a comparative analysis of the relative error
between the new and existing solutions. The digital data reveal a remarkable agreement
between these results, especially with those in \cite{liu2023semi} for one-sided FPDE problems.
The experimental validation demonstrates satisfactorily the reliability and accuracy of
new finite difference scheme (\ref{e2-9}) for predicting quenching time. 

\begin{table}
\begin{center}
 \tabcolsep 0.1in\small
\caption{A comparison between the quenching time $T_{a}$ and $T_{a}^{S},T_{a}^{M},T_{a}^{L}.$ ($\theta=1;\,a=\pi$).}
\vspace{2mm}
\begin{tabular}{lllllll}
 \hline
  $b$ &  $-2.0$ & $-0.9$ & $-0.4$ & $0.4$ & $0.9$ & $0.95$\\\hline
  $T_{a}$& 0.5846 & 0.5577 & 0.5465  & 0.5300 & 0.5216 & 0.5209 \\
  $T_{a}^{S}$ &0.5680 &0.5490 &0.5420 & 0.5420 & 0.5320 & 0.5320\\
  $T_{a}^{M}$ &0.5880 &0.5590 &0.5470 & 0.5280 & 0.5110 & 0.5080\\
  $T_{a}^{L}$& 0.5850 & 0.5580 & 0.5468  & 0.5304 & 0.5220 & 0.5212 \\[2pt]\hline
 \end{tabular}
 \label{tab3}
 \end{center}
\end{table}

\begin{table}
\begin{center}
 \tabcolsep 0.05in\small
\caption{A comparison of the quenching time $T_{a}$ and $T_{a}^{S},T_{a}^{M},T_{a}^{L}.$ ($\theta=1;\,a=2$). }

\vspace{2mm}
\begin{tabular}{llllllllllll}
 \hline
  $b$ & $-0.9$ &$-0.8$ &$-0.6$ &$-0.4$ &$-0.2$ &$0.2$ & $0.4$ &$0.6$ & $0.8$& $0.9$  \\\hline 
  $T_{a}$ & 1.1517 & 1.0836 &0.9765 &0.8956 &0.8315 &0.7365 &0.7013 &0.6728 &0.6506 &0.6417  \\       
  $T_{a}^{S}$ & 1.0870 & 0.9870 &0.9160 &0.8600 &0.8150 &0.7500 &0.7300 &0.7130 &0.7060 &0.7060  \\
  $T_{a}^{M}$ & 1.1600 &1.0900 &0.9810 &0.8990 &0.8340 &0.7320 &0.6880 &0.6470 &0.6030&0.5760 \\
  $T_{a}^{L}$ & 1.1522 & 1.0840 &0.9770 &0.8962&0.8320 &0.7370 &0.7018 &0.6732 &0.6510 &0.6420 \\[2pt]\hline
 \end{tabular}
 \label{tab4}
 \end{center}
\end{table}

\begin{table}
\begin{center} \tabcolsep 0.1in\small
\caption{A data observation of the monotonic convergence of
$T_{a}$ as $a\rightarrow\infty$  ($b=-0.4, \sigma=2,\theta=1$).}
\vspace{2mm}

\begin{tabular}{ll   ll   ll   ll}
 \hline
 $a$  &  $T_{a}$  &  $a$  &  $T_{a}$  &  $a$  &  $T_{a}$  &  $a$  &  $T_{a}$ \\\hline
 1.7 & $2.4888$ &  2.0& 0.8956  & 4.0 & 0.5136  & 10.0 & 0.5010  \\
 1.8 & $1.3656$ & 3.0 & 0.5586  &5.0& 0.5037  & 20.0 & 0.5010 \\
 1.9 & 1.0521  & $\pi$ & 0.5472  & 8.0 &0.5010 & 40.0 & 0.5010  \\[2pt]\hline
 \end{tabular}
 \label{tab5}
 \end{center}
\end{table}

\begin{figure}
\begin{center}
\epsfig{file=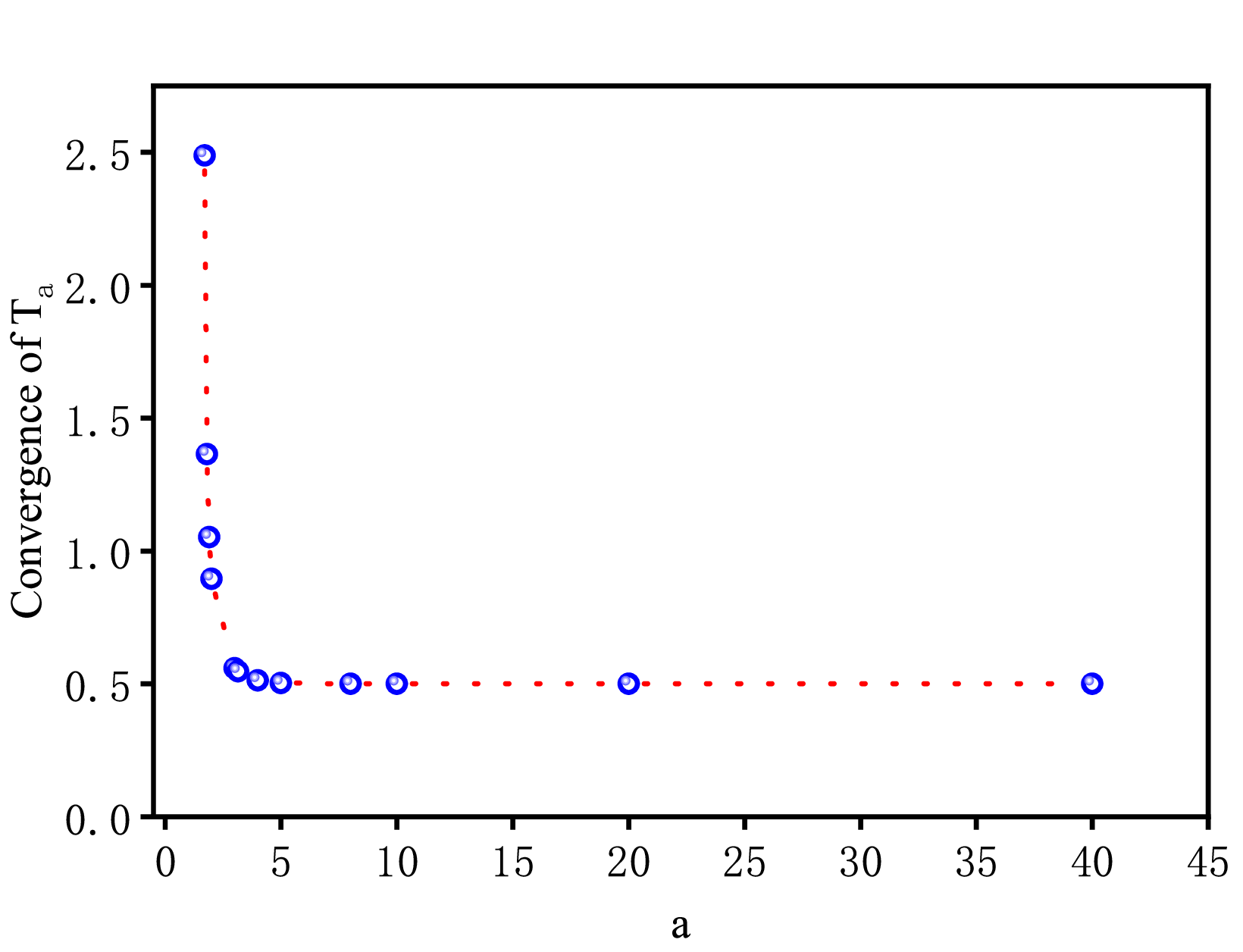,width=7.5cm,height=5.7cm}
\caption{A monotonic convergence of $T_{a}$ as $a\rightarrow\infty.$ }
\label{fig4}
\end{center}
\end{figure}

Furthermore, Table \ref{tab5} indicates a monotonic convergence of the quenching time $T_{a}.$
It is found that $T_{a}$ monotonically decreases and converges to 0.5 as
the interval size $a\rightarrow\infty.$
Our experiments also indicate that, as the size of the spacial interval increases,
$T_{a}$ gradually converges towards a fixed value of $0.5.$ 
The finding aligns with numerical results outlined in \cite{zhu2020note}.
Additionally, Fig. \ref{fig3} offers a simulation of the relative errors of the ratios between our quenching time and $b$
compared to the well-accepted known results. They seem to be highly satisfactory.
Again, the monotone convergence phenomenon is shown in Fig. \ref{fig4}.

\begin{table}
\begin{center}
 \tabcolsep 0.04in\small
\caption{The variation of quenching time $T_{a}$ with respect to the parameter $b$ ($\theta=1$).}

\vspace{2mm}
\begin{tabular}{llllllllllll}
 \hline
  $b$ & $-1.0$ &$-0.6$ &$-0.4$ &$0$ &$0.4$ &$0.6$ & $1.0$ &$1.4$ & $1.8$& $2.0$  \\\hline 
  $T_{2}(\sigma=1.8)$  &0.8772 &0.7740 &0.7341 &0.6669 &0.6189 &0.6042 &0.5907 &0.5925 &0.6060 &0.6171 \\
 $T_{2}(\sigma=2)$  &1.2360 &0.9777 &0.8967 &0.7803 & 0.7020 & 0.6735 &0.6351 &0.6183 &0.6219 &0.6309 \\
 $T_{\pi}(\sigma=1.8)$ &0.5484 &0.5421 &0.5388 &0.5319 &0.5262 &0.5241 &0.5211 &0.5196 &0.5190 &0.5190 \\
 $T_{\pi}(\sigma=2)$  &0.5607 &0.5517 &0.5472 &0.5385 &0.5307 &0.5271 &0.5208 &0.5160 &0.5130 &0.5118 \\ [2pt]\hline
 \end{tabular}
 \label{tab6}
 \end{center}
\end{table}

\begin{table}
\begin{center}
 \tabcolsep 0.06in\small
\caption{Quenching time $T_{a}$ versus fractional order of derivative $\sigma$ ($b=-0.4,\theta=1$).}

\vspace{2mm}
\begin{tabular}{llllllllllll}
 \hline
  $\sigma$ & $1.65$ &$1.7$ &$1.75$ &$1.8$ &$1.85$ &$1.9$ & $1.95$ &$2.0$ \\\hline 
   $T_{1.8}$  &0.7230 &0.7698 &0.8253 &0.8916 &0.9723 &1.0725 &1.1994 &1.3656 \\
 $T_{2}$  &0.6477 &0.6735 &0.7023 &0.7341 &0.7695 &0.8082 &0.8508 &0.8967 \\
 $T_{\pi}$ &0.5316 &0.5340 &0.5364 &0.5388 &0.5412 &0.5433 &0.5454 &0.5472 \\
 $T_{10}$  &0.5028 &0.5025 &0.5022 &0.5019 &0.5016 &0.5013 &0.5010 &0.5010 \\ [2pt]\hline
 \end{tabular}
 \label{tab7}
 \end{center}
\end{table}

\begin{figure}
\begin{center}
\epsfig{file=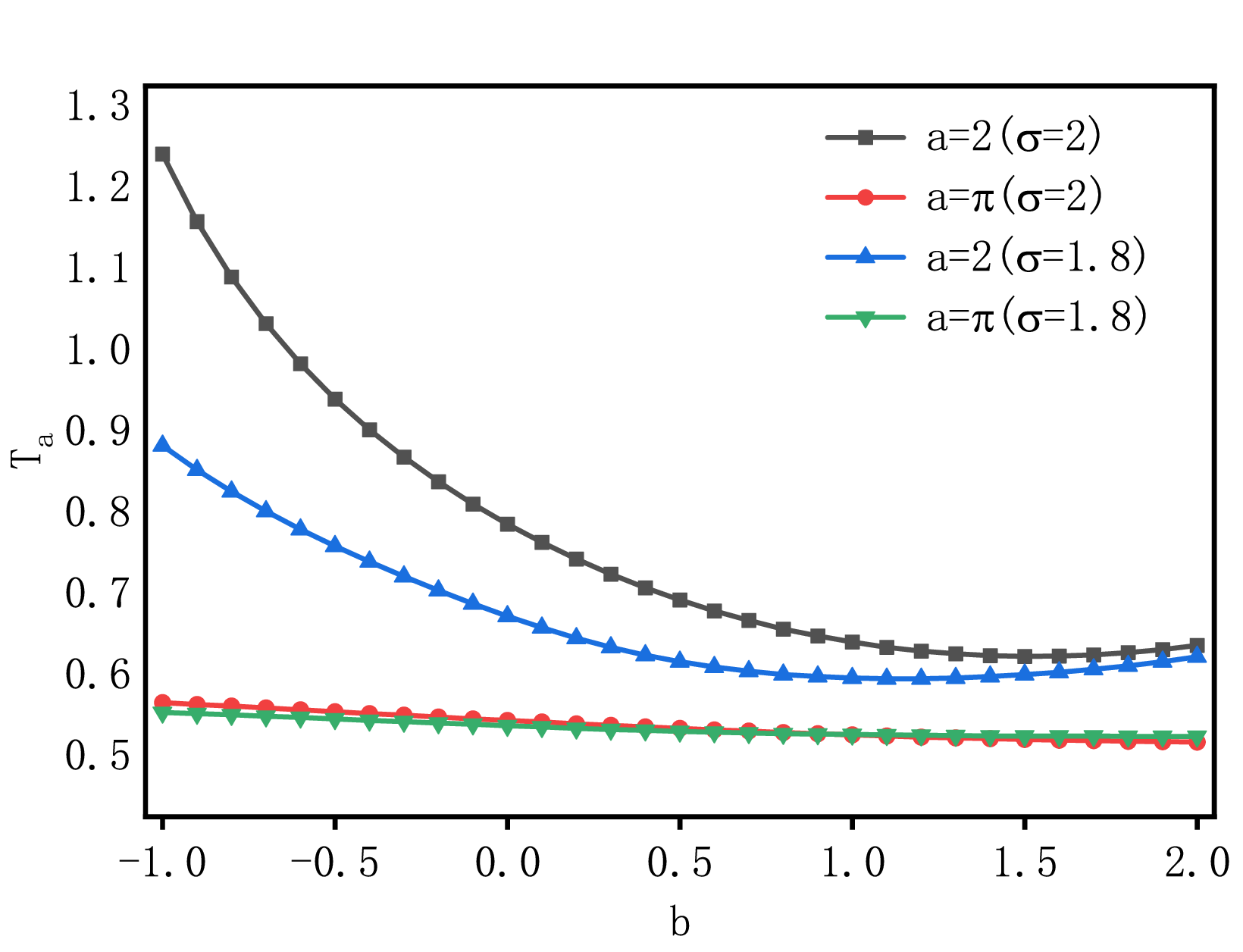,width=7.0cm,height=5.4cm}
\epsfig{file=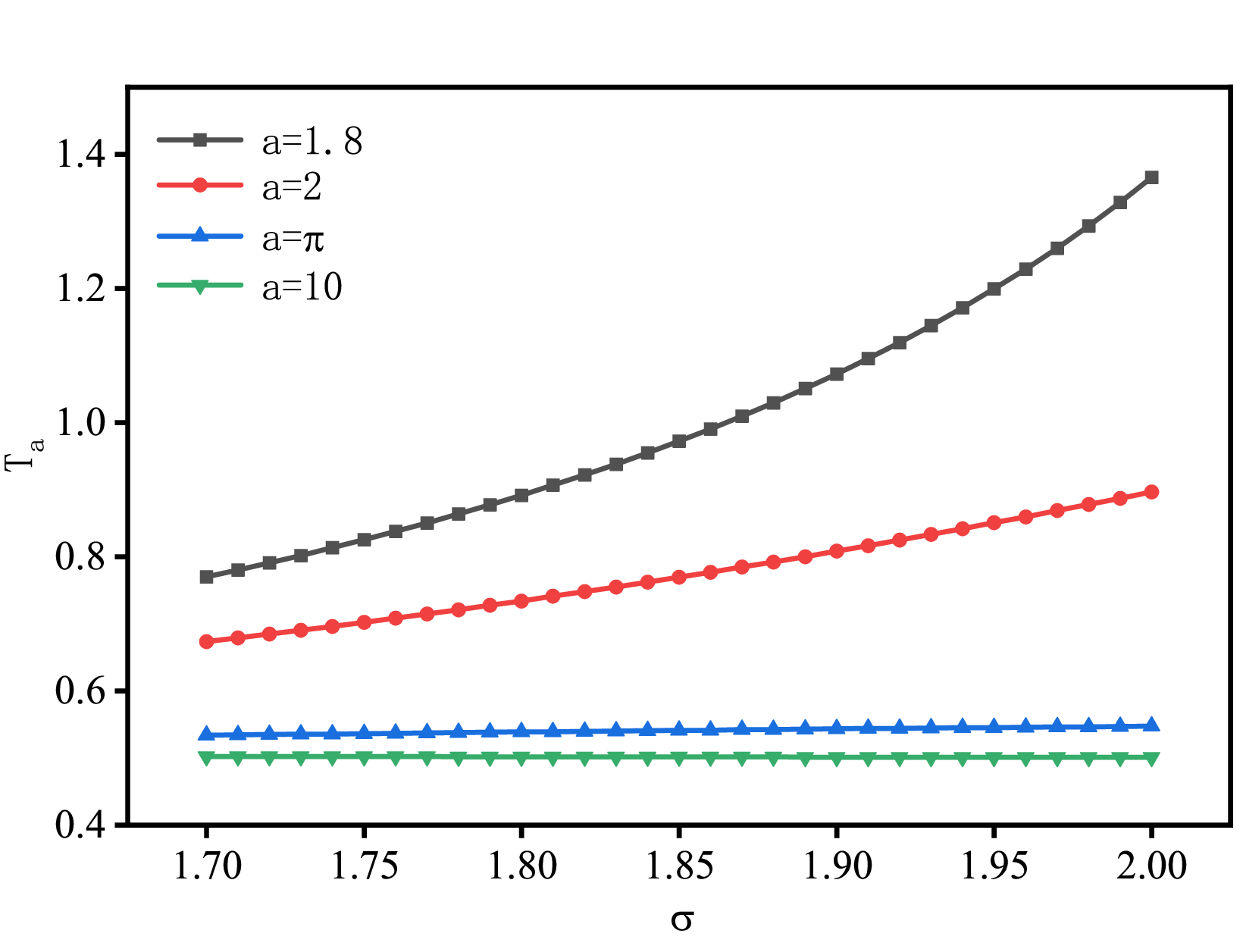,width=7.0cm,height=5.4cm}
\caption{LEFT: Dependence of $T_{a}$ on $b.$ RIGHT: Dependence of $T_{a}$ on the fractional derivative order $\sigma.$}
\label{fig5}
\end{center}
\end{figure}

Table \ref{tab6} outlines correlations between $T_{a}$ and the coefficient $b$ when
different fractional orders are used. In Fig. \ref{fig5} (left), it becomes evident that $T_{a}$ monotonically decreases
as $b$ increases, and it ultimately converges to $0.5.$ Noticeably, the quenching time $T_{a}$ as $a=2$ is significantly
greater than that when $a=\pi.$

Furthermore, Table \ref{tab7} presents results from the numerical analysis
about the quenching time $T_{a}$ in connection with the fractional order $\sigma.$
As depicted in Fig. \ref{fig5} (right), it is clearly visible that $T_{a}$ decreases when the
value of $\sigma$ is reduced.
The phenomenon precisely repeats that observed in \cite{liu2023semi}, while the
current paper offers a much clearer profile how the value of $T_{a}$ changes as the fractional order
varies. The experiment provides a valuable insight into the richer characteristics of fractional order
quenching modeling problems.


\subsection{Simulation experiment C: quenching location of the solutions}
\label{sec: sec4-3}

Let us again consider the semi-adaptive scheme (\ref{e2-9}) while new parameters $\sigma=2, a=\pi$
and $\theta=1.0$ are used.
The quenching times are found to be $0.5490$ and $0.5284$ when $b = -0.5,~0.5,$ respectively.
Furthermore, the corresponding quenching position approximations are detected at $x^{*}=1.7279$ and $x^{*}=1.3823,$
respectively. Noticeably, when $\sigma=2,$ that is, the order is an integer, our computed results are
highly agree with known values \cite{sk2001,sheng1999compound,Josh3,mooney1996implicit}.

\begin{figure}
\begin{center}
\epsfig{file=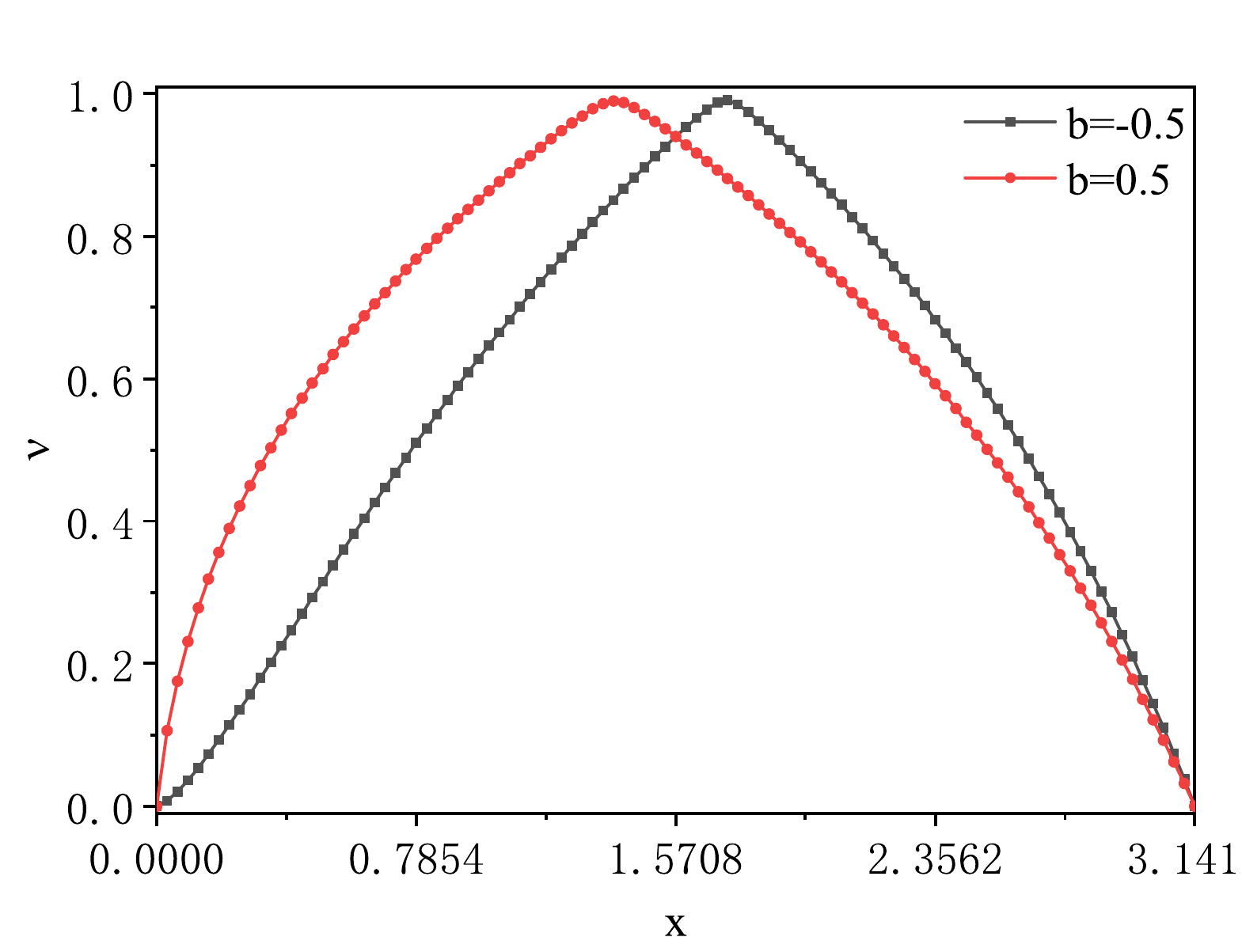,width=7.3cm,height=5.2cm}
\epsfig{file=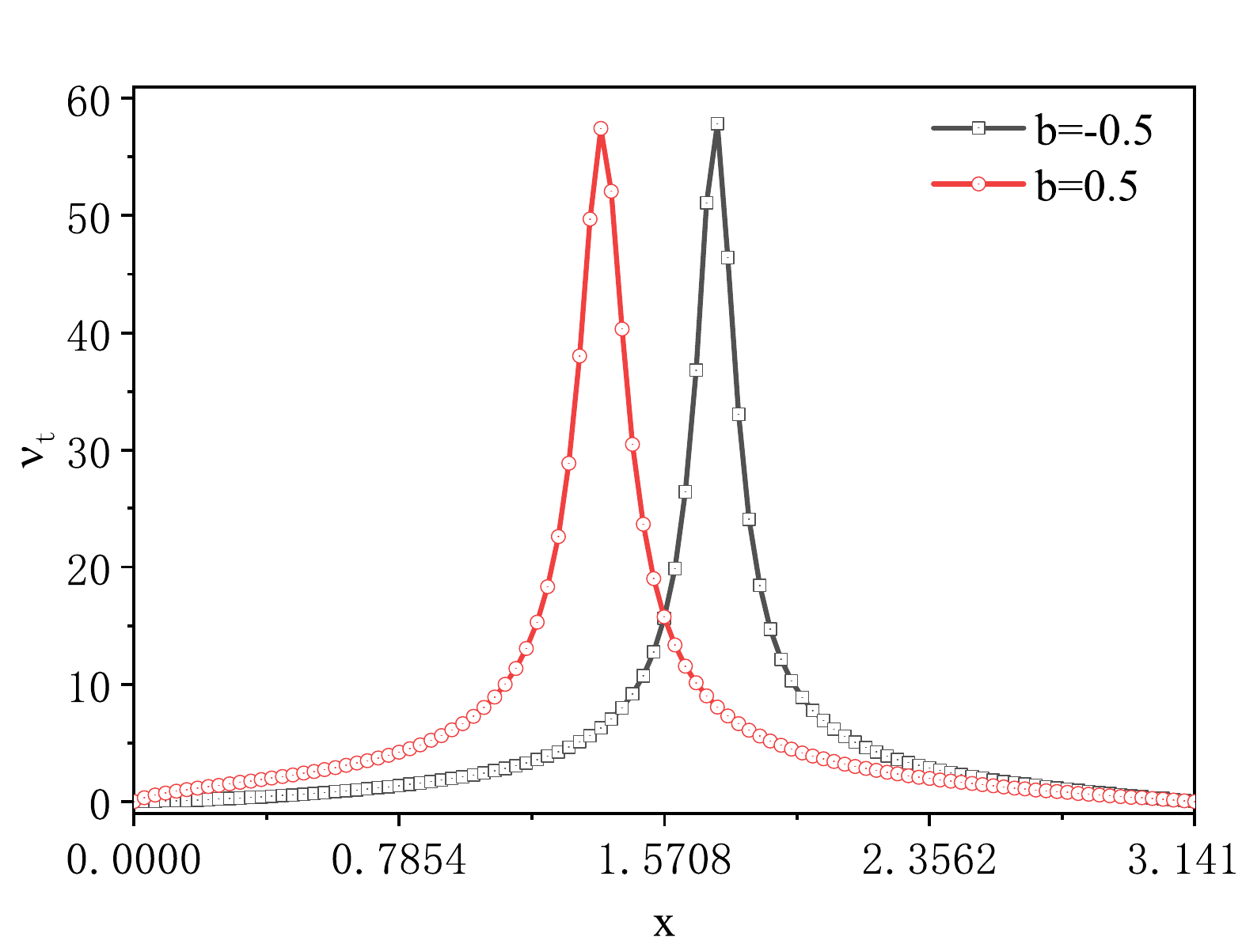,width=7.3cm,height=5.2cm}
\epsfig{file=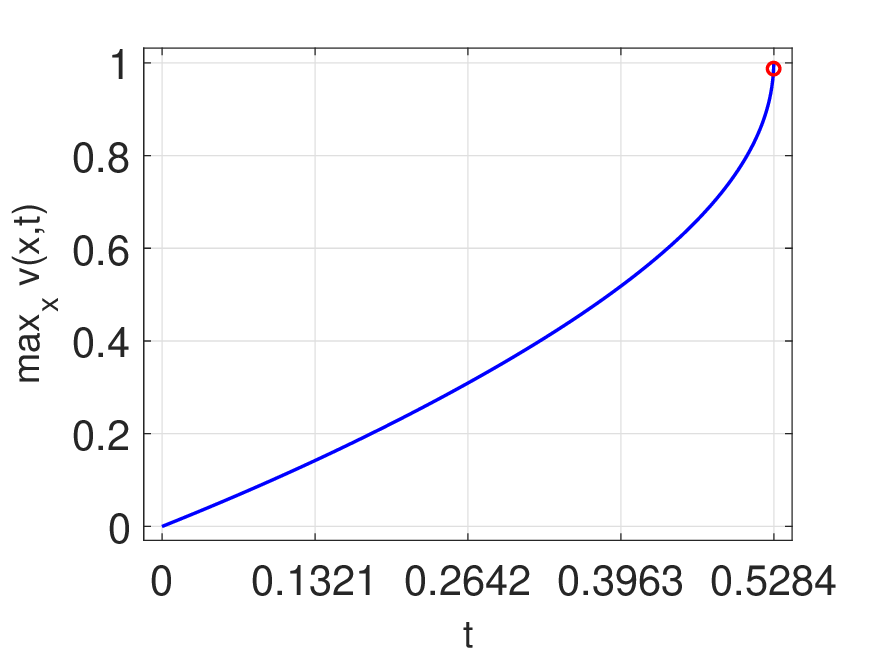,width=7.4cm,height=5.2cm}
\epsfig{file=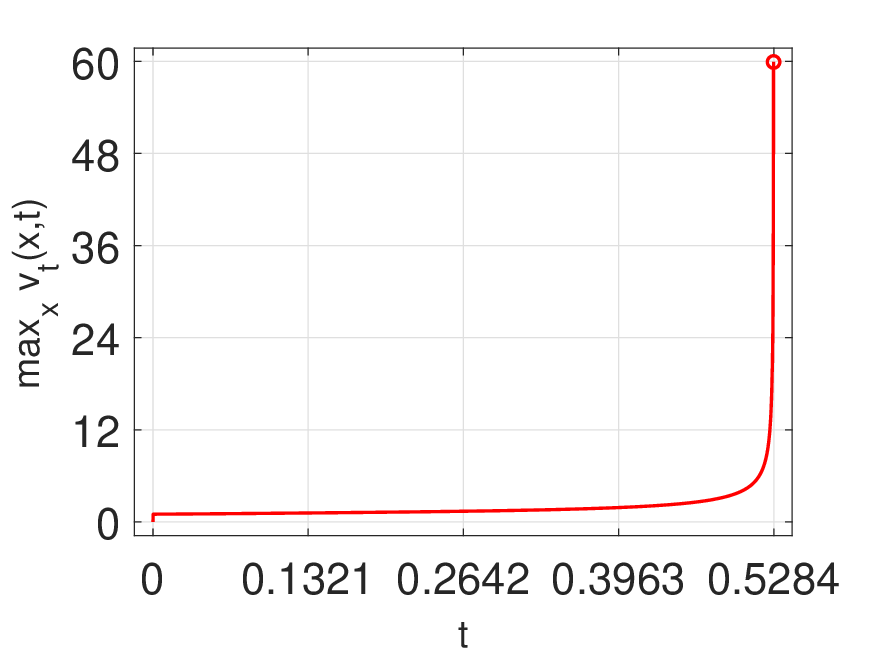,width=7.4cm,height=5.2cm}
\caption{TOP: Cross sections of $v$ and $v_{t}$ immediately prior to quenching.
BOTTOM: Trajectories of the peak values of $v$ and $v_t$ before quenching.
Note that quenching occurs at time $T_{\pi}\approx 0.5284~(\sigma=2,a=\pi,~b=0.5).$ }
\label{fig6}
\end{center}
\end{figure}

On the other hand, Fig. \ref{fig6} (top) shows distributions of the intersections of the numerical solution $v$
and its temporal $v_{t}$ immediately prior to quenching. The quenching singularity occurs slightly earlier if $b=0.5$ as
compared with that of $b=-0.5.$ In other words, the larger the value of $b$ is, the earlier the quenching may appear.
Furthermore, due to a possible influence of the variable convection coefficient $c(x,t),$ the numerical solution $v$ and
the corresponding temporal derivative $v_{t}$ exhibit a non-axisymmetric behavior, while reaching their extrema
approximately at the spacial location $x=1.5708.$

Fig. \ref{fig6} (bottom) records evolution profiles of the extrema of $v$ and $v_{t}$ from $t=0$ up to the
quenching time $T_{\pi}.$
It is clear that both trajectories increase smoothly  and monotonically till the quenching blow-up is reached, that is,
as $t\approx T_{\pi}\approx 0.5284.$ At a quenching singularity, the maximal value of $v$
approaches quickly the unit, whereas the superior of $v_{t}$ surges exponential and grows ultimately to
59.9161 in our experiment. The phenomena are consistent with those reported in the latest
nonlinear Kawarada equation studies (see \cite{sheng2023nonconventional,padgett2018quenching} and
references therein).

\begin{figure}
\begin{center}
\epsfig{file=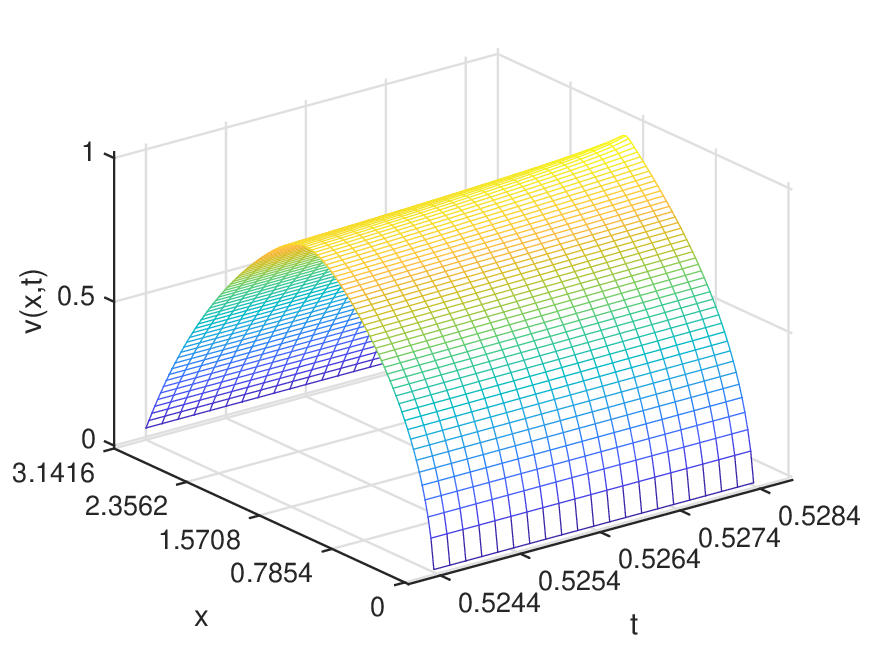,width=7.3cm,height=5.2cm}
\epsfig{file=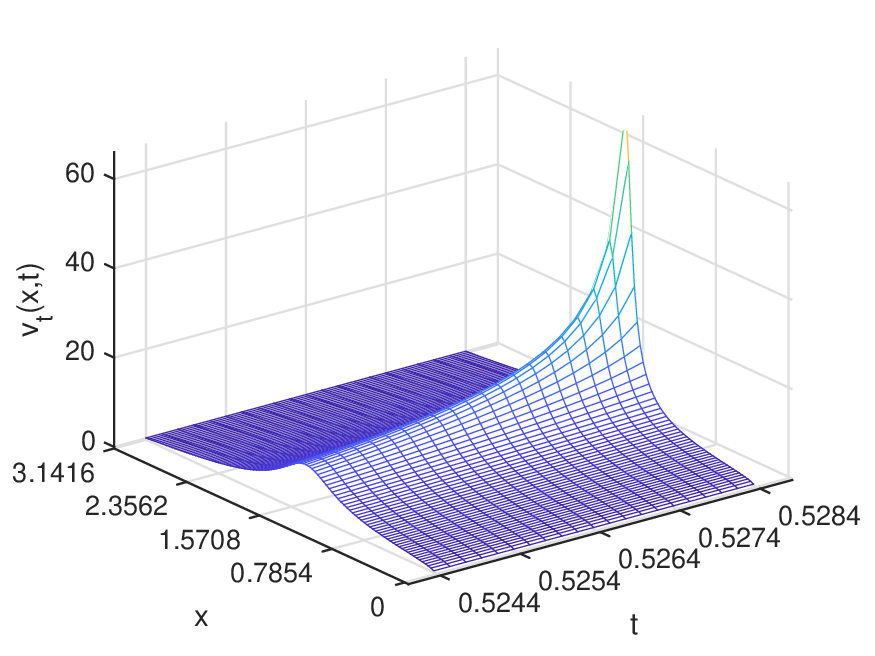,width=7.3cm,height=5.2cm}
\caption{3D surface plots of $v(x,t)$ and $v_{t}(x,t)$ on the 20 temporal steps immediately
before quenching. While $v$ peaks at 1, $v_{t}$ peaks at 59.9161, approximately,
($\sigma=2,a=\pi,b=0.5$).}
\label{fig7}
\end{center}
\end{figure}

To see more precisely simulations of quenching solutions and to examine the reliability of the algorithm,
we now choose a case where $\sigma=2,a=\pi,b=0.5$ are employed. Within the temporal interval
$t\in [0.5244, 0.5284]$ immediately prior to quenching, we generate three-dimensional
surface plots of $v$ and $v_{t},$ respectively. Apparently, the quenching location is approximately $x^{*}=1.3823.$
As the quenching moment is getting close, the maximal value of the numerical solution $v(x,t)$ approaches to one,
while the superior of $v(x,t)$ rapidly surges to $59.9161$ in our experiments. Fig. \ref{fig6} is devoted to a planar
cross-section for Fig. \ref{fig7}. They both offer accurate and reliable visualizations of the solutions and
their evolutionary profiles. The explorations allow us to gain decent insights into quenching solutions' rich
dynamic behavior.

\begin{table}
\begin{center}
 \tabcolsep 0.12in\small
\caption{Peak values of $v(x,t)$ and $v_{t}(x,t)$ at various quenching locations $x^*$
($a=2,\theta=1.0,b=-0.5$). }

\vspace{2mm}
\begin{tabular}{llllllllllll}
 \hline
  $\sigma$ & $v^{*}(x,t)$ & $v^{*}_{t}(x,t)$ & $T_{a}$ & $x^{*}$ \\\hline
  $1.7$ & 0.9930 & 52.2753 & 0.6862 & 1.20 \\
  $1.75$ & 0.9920 & 47.3515 & 0.7178 & 1.20 \\
  $1.8$ & 0.9987 & 52.8629 & 0.7530 & 1.18 \\
  $1.85$ & 0.9914 & 51.6327 & 0.7920 & 1.18 \\
  $1.9$ & 0.9965 & 43.4455 & 0.8350 & 1.16 \\[2pt]\hline
 \end{tabular}
 \label{tab8}
 \end{center}
\end{table}

\begin{figure}
\begin{center}
\epsfig{file=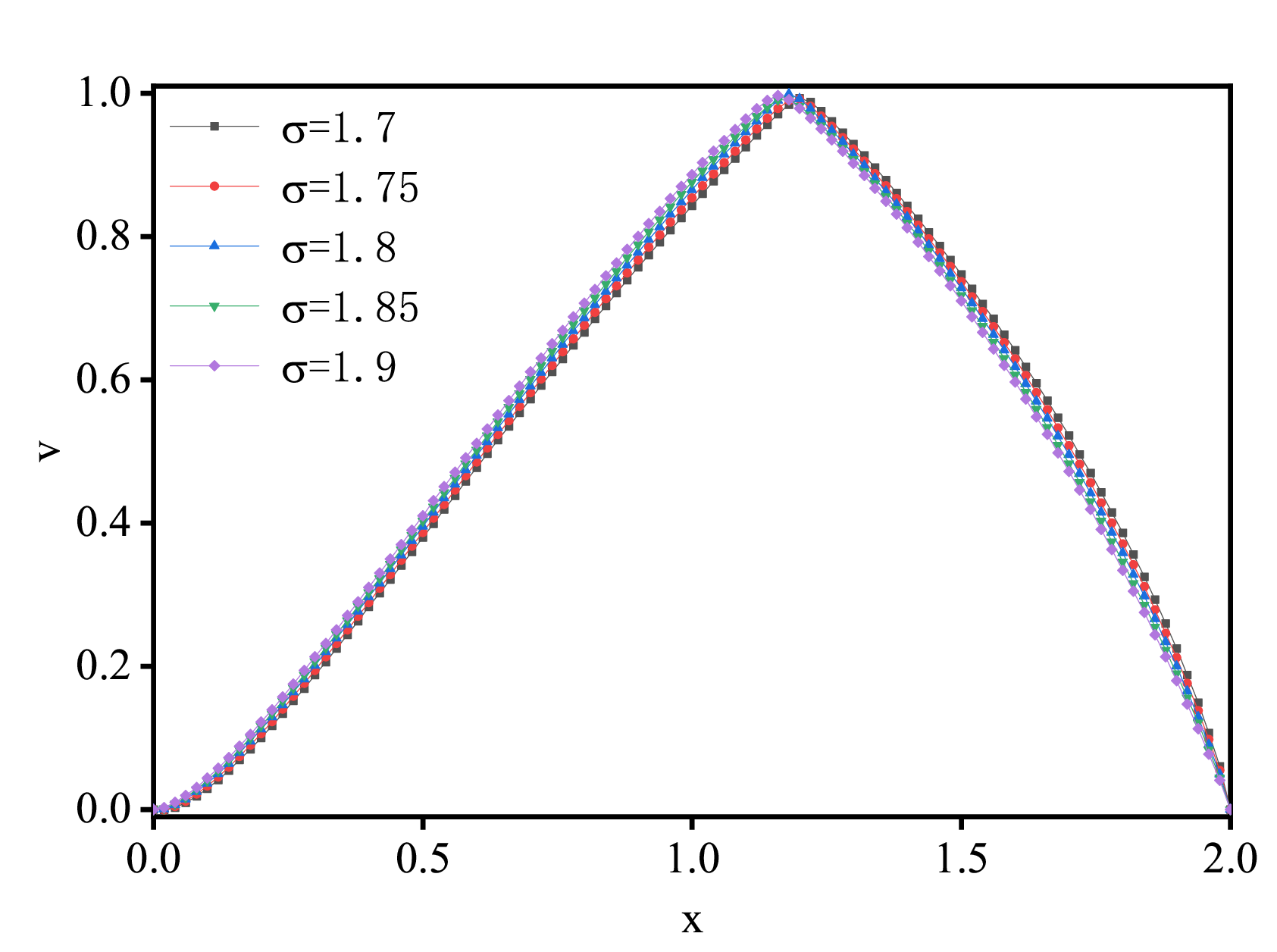,width=7.1cm,height=5.5cm}
\epsfig{file=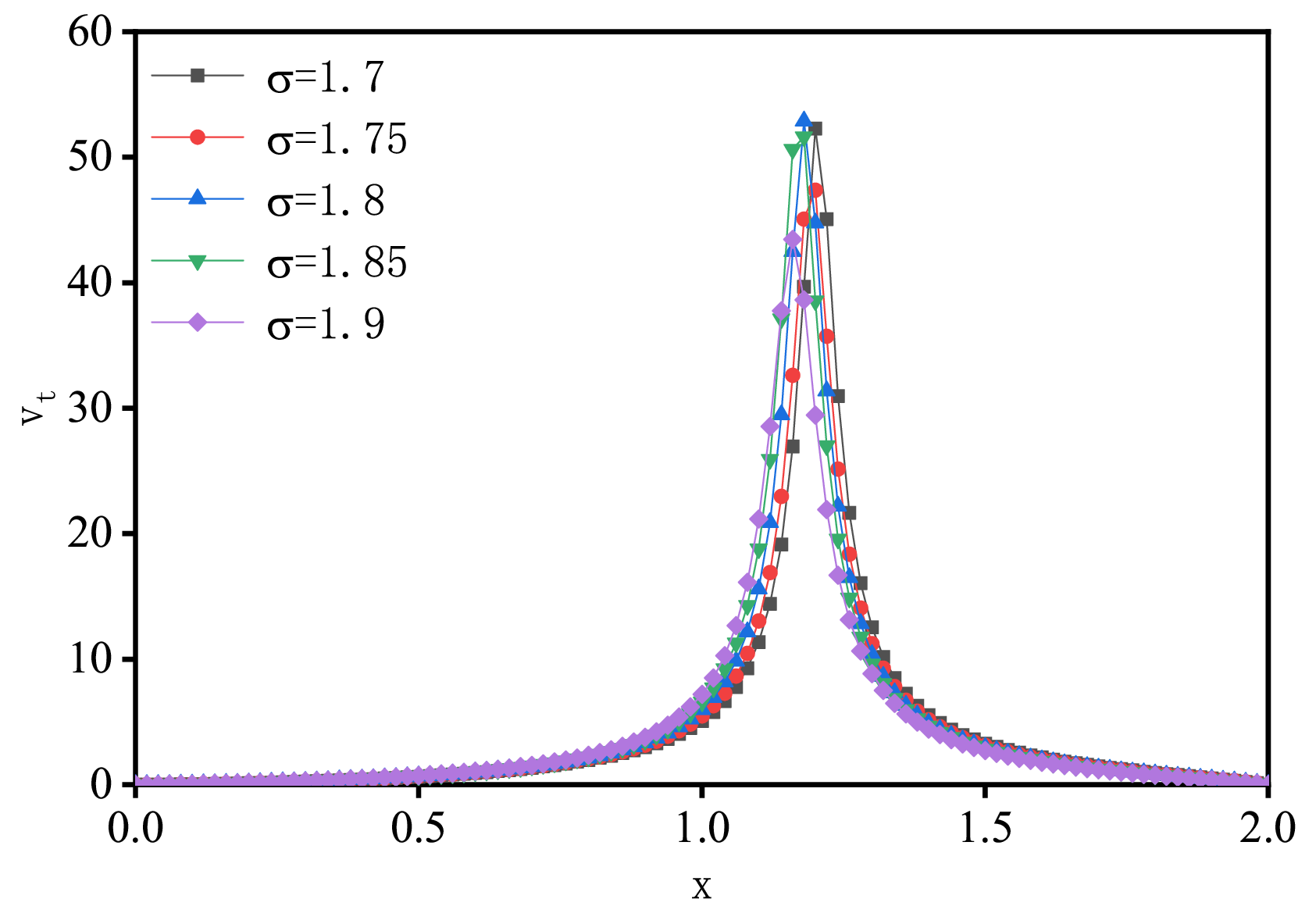,width=7.1cm,height=5.2cm}
\caption{Cross sections of  $v$ and $v_{t},$ respectively,
immediately before quenching ($a=2,\,b=-0.5$).}
\label{fig8}
\end{center}
\end{figure}

Now, for $a=2, \theta=1.0,$ and a negative $b=-0.5,$ we list peak values of the solution functions $v$ and $v_{t}$
at quenching time $T_{a}$ and quenching location $x^{\star}$ in Table \ref{tab8}.
Five distinct non-integer values of $\sigma$ ranging from 1.7 to 1.9
are used. For the same parameters used in Table \ref{tab8}, Fig. \ref{fig8} draws five cross section curves for
each of $v$ and $v_{t}$ with respect to the spatial variable $x$ immediately before the quenching
corresponding to aforementioned five $\sigma$ values.
It can be found clearly that as the fractional order increases from 1.7 to 1.9, the quenching location $x^*$ shifts to
the right slowly. It is noticed that at $\sigma=1.8,$ we have approximately $x^{*}\approx 1.18,$ where the rate of
change function \cite{sk2001} $v_{t}$ attains its maximum 52.8629.

\begin{figure}
\begin{center}
\epsfig{file=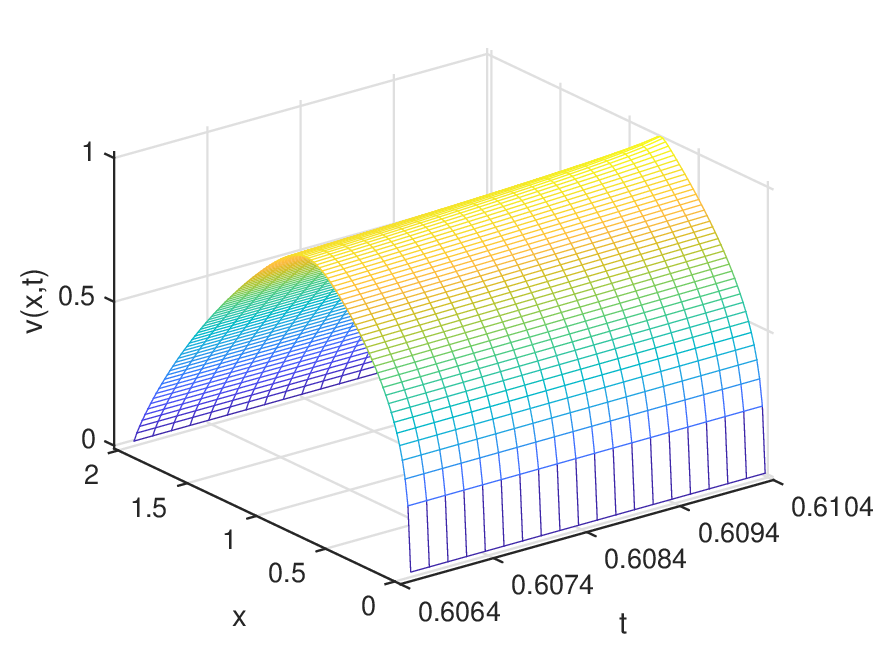,width=7.3cm,height=5.2cm}
\epsfig{file=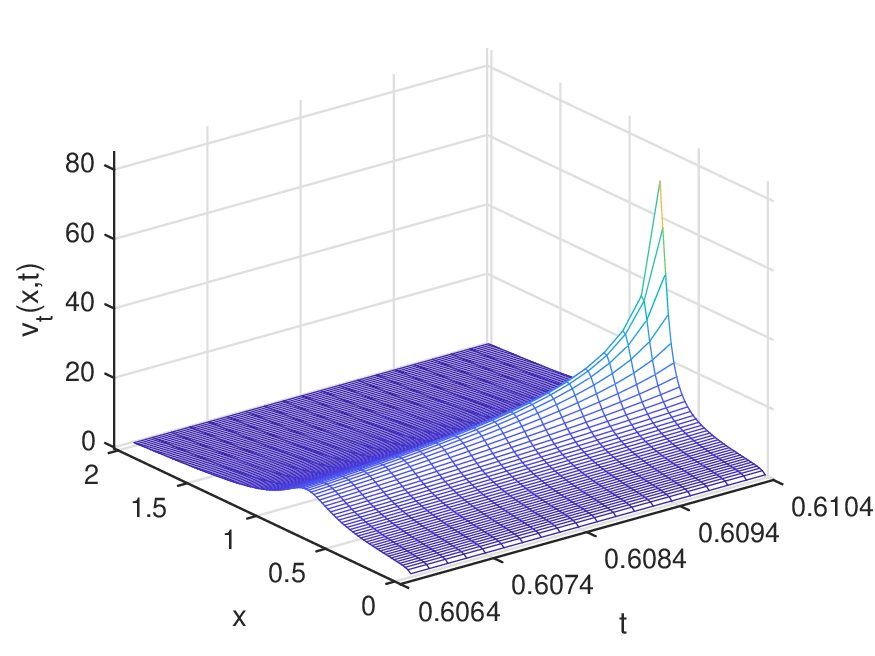,width=7.3cm,height=5.2cm}
\epsfig{file=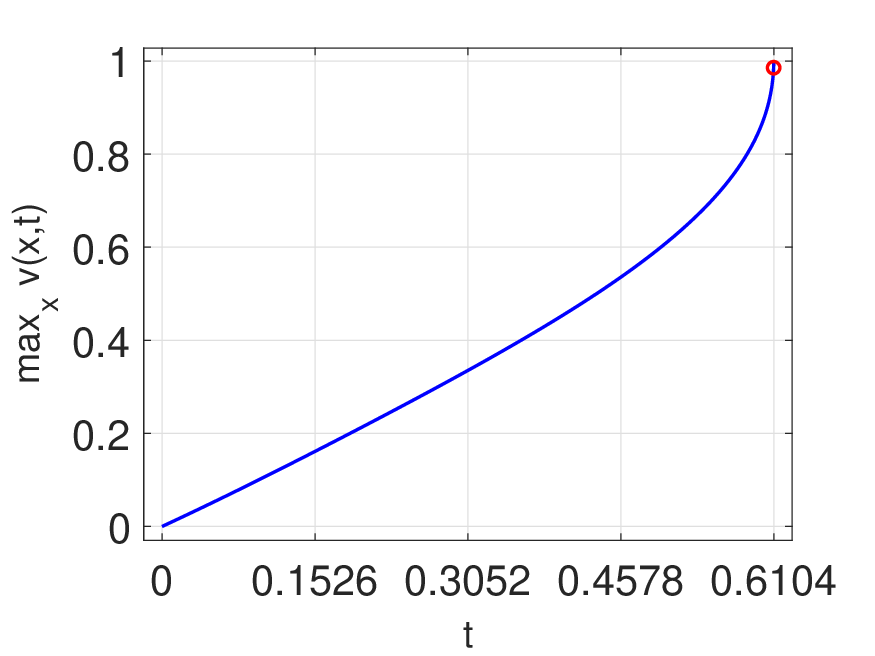,width=7.3cm,height=5.2cm}
\epsfig{file=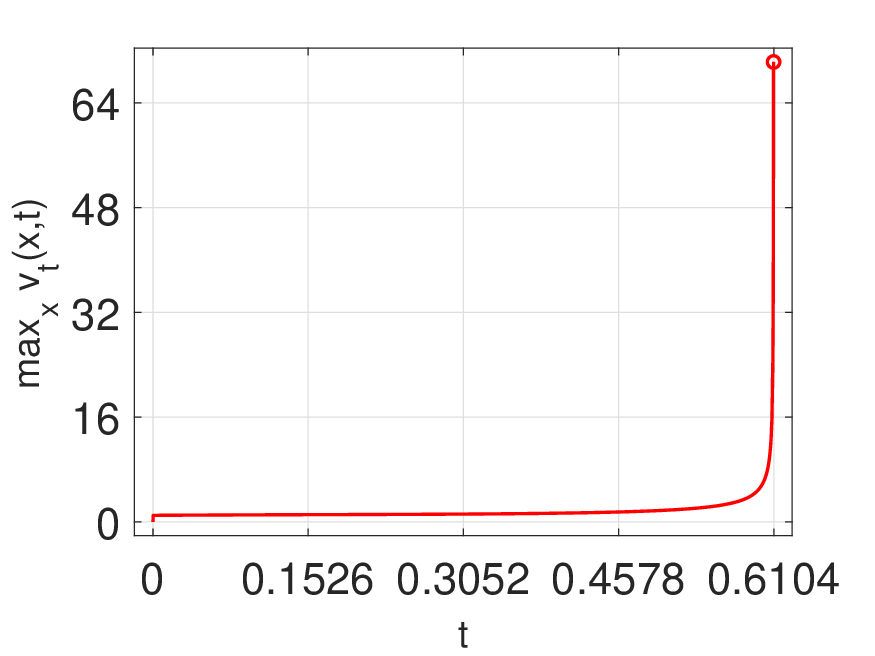,width=7.3cm,height=5.2cm}
\caption{TOP: Surfaces of $v(x,t)$ and $v_{t}(x,t)$ immediately before quenching.
BOTTOM: Trajectories of peak values of $v$ and $v_t$ before quenching. Note that
quenching occurs at time $T_{2}\approx 0.6104~ (\sigma=1.8,a=2,b=0.5,\theta=1.0).$}
\label{fig9}
\end{center}
\end{figure}

Moreover, the quenching time obtained is $T_{2}=0.6104$ when parameters $\sigma=1.8,a=2,b=0.5$ are used.
The 3D profiles of the numerical solution $v$ and corresponding temporal derivative $v_{t}$ with a short
time interval prior to quenching are given in Fig. \ref{fig9}.
We may observe that $\max\limits_{0\leq x\leq a} v(x,t)$ converges quickly to the unity as $t\rightarrow 0.6104,$ while
$\sup\limits_{0< x< a} v_{t}(x,t)$ approaches 70.2520. exponentially. The phenomenon may suggest that the
value of $\sup\limits_{0< x< a} v_{t}(x,t)$ becomes unbounded at the quenching.
Fig. \ref{fig9} (bottom) is devoted to  profiles of the maximal values of $v$ and $v_{t}$ within the last few temporal steps
before quenching ($0.6064\leq t\leq 0.6104$).  The surface plots also highlight the sensitivity of the quenching
phenomena associated with fractional order quenching problems. It is worth to point out that profiles of $v$ and
$v_{t}$ shown in our simulation experiments again exhibit high degrees of agreements with
data given and analyzed in multiple recent publications (see \cite{sk2001,padgett2018quenching,liu2023semi}
and references therein, for instance).

\begin{table}
\begin{center}
 \tabcolsep 0.12in\small
\caption{Peak values of $v(x,T_a)$ and $v_{t}(x,T_a)$ at quenching location $x^*$ ($a=2,\,\sigma=1.8$).}

\vspace{2mm}
\begin{tabular}{llllllllllll}
 \hline
  $b$ & $v^{*}(x,t)$ & $v^{*}_{t}(x,t)$ & $T_{a}$ & $x^{*}$ \\\hline
  $-1.0$ & 0.9921 & 59.2078 & 0.8765 & 1.32 \\
  $-0.6$ & 0.9988 & 51.2375 & 0.7736 & 1.22 \\
  $-0.2$ & 0.9904 & 46.5946 & 0.6982 & 1.08 \\
  $0$ & 0.9912 & 51.6249 & 0.6666 & 1.00 \\
  $0.2$ & 0.9908 & 43.0671 & 0.6396 & 0.90 \\
  $0.6$ & 0.9922 & 39.8836 & 0.6038 & 0.70 \\
  $1.0$ & 0.9985 & 58.3906 & 0.5904 & 0.52 \\[2pt]\hline
 \end{tabular}
 \label{tab9}
 \end{center}
\end{table}

\begin{figure}
\begin{center}
\epsfig{file=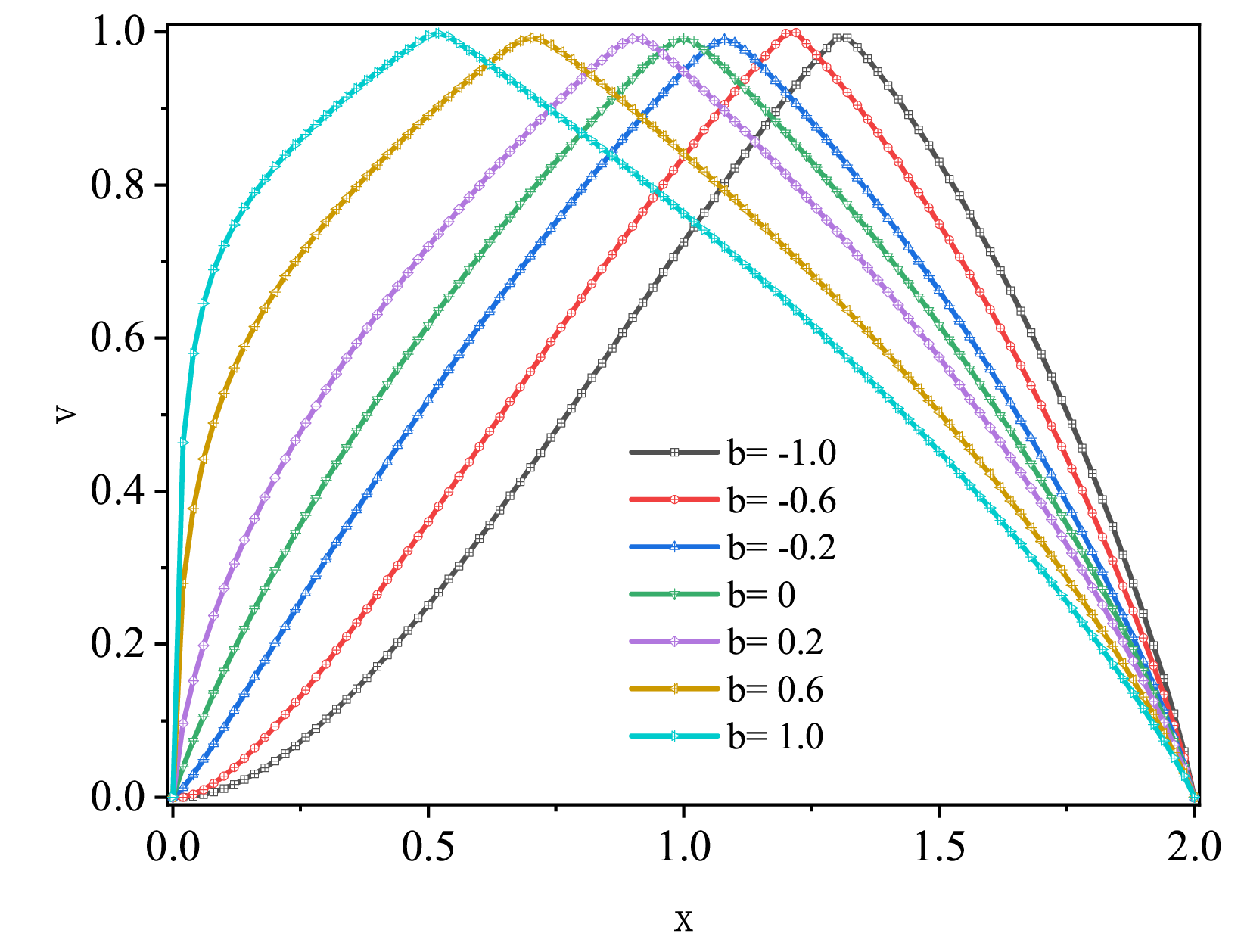,width=7.3cm,height=5.0cm}
\epsfig{file=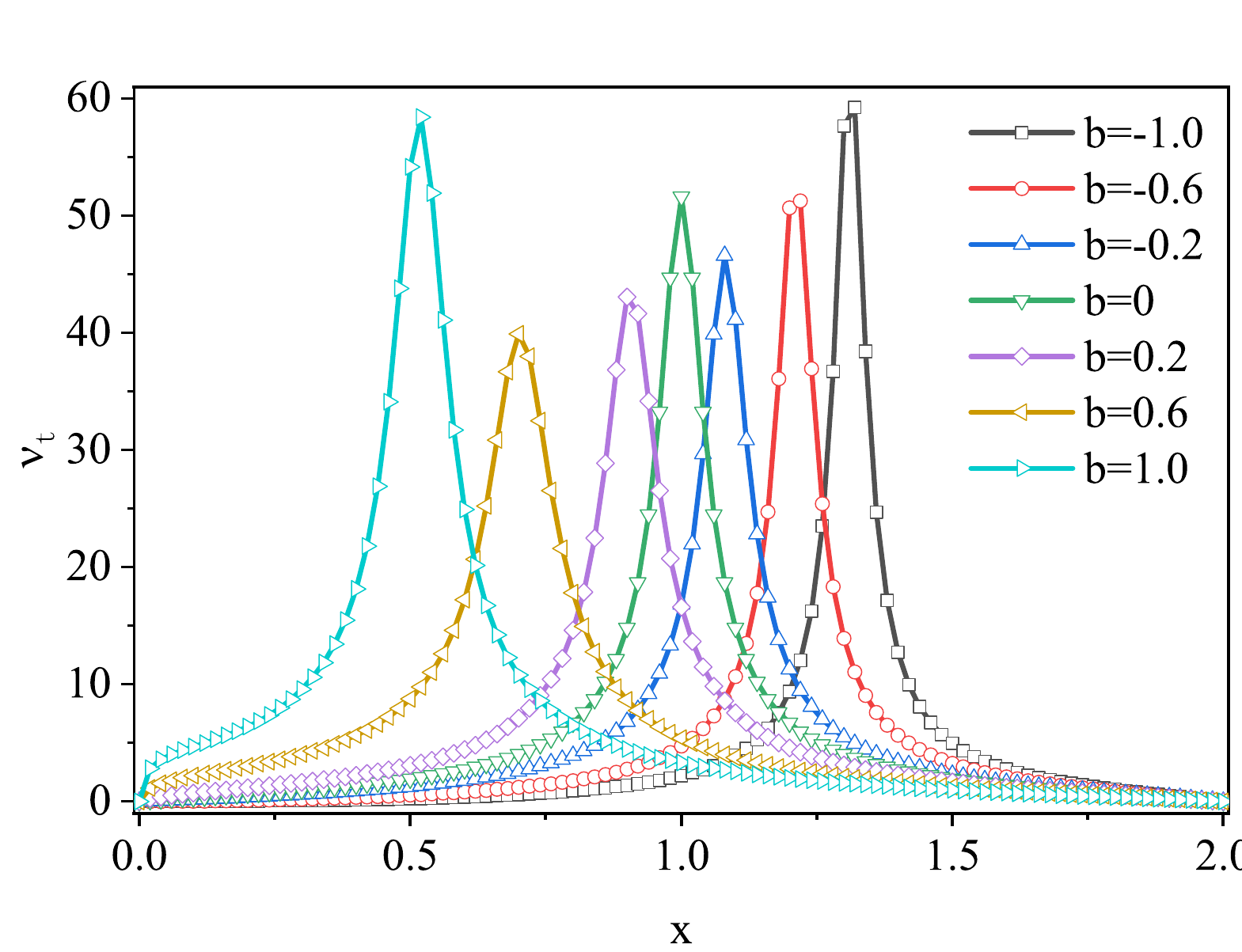,width=7.3cm,height=5.4cm}
\caption{Shifts of quenching locations with respect to the parameter $b$ ($a=2,\sigma=1.8$).}
\label{fig10}
\end{center}
\end{figure}

Advanced with $\sigma=1.8, a=2, \theta=1.0,$ and variable $b$ values from -1.0 to 1.0,
Table \ref{tab9} summarizes maximal values of the numerical solution $v(x,t)$ and
rate of change function $v_{t}(x,t)$ priori to the quenching time $T_{a}.$
Quenching positions $x^{*}$ corresponding to different values of $b$ are detected.
Fig. \ref{fig10} continuously illustrates visible shifts of quenching positions as the diffusion coefficient $b$ changes.
In fact, we find that the quenching location decreases monotonically as $b$ increases within the
above range. Geometrically speaking, the quenching location gradually shifting from the right
towards the left in the spacial interval $(0,a).$ We are particularly interested in the phenomenon that
when $b$ reaches the value of -1.0, the peak value of $v_{t}$ approaches $59.2078$ rapidly at $x^*.$
This implies the sensitivity of a rational order quenching solution with respect to $b.$
We may also conclude that both $v(x,t)$ and $v_{t}(x,t)$ of the two-sided fractional convection-diffusion
quenching model (\ref{e2-1}) may not be spatially symmetric. This confirms several recent discoveries
reported for nonlocal quenching problems \cite{zhu2023simulation,liu2023semi}.

\subsection{Simulation experiment D: orders of convergence of the solutions}
\label{sec: sec4-4}

To study the convergence of the semi-adaptive method (\ref{e2-9}), our numerical experiments are extended for
estimating the order of convergence the numerical solution $v(x,t)$ and its temporal derivative $v_{t}(x,t)$
obtained from (\ref{e2-9}).
We set the spatial step $h=a/100$ and and temporal  step $\tau=10^{-4}$
and execute $(\ref{e2-9})$ for numerical solutions $v_{h}^{\tau}, v_{{h}/{2}}^{\tau}, v_{{h}/{4}}^{\tau},$ respectively.
Consider the solutions within the final 20 temporal steps before quenching, and
evaluate corresponding derivatives $(v_t)_{h}^{\tau}, (v_t)_{{h}/{2}}^{\tau}, (v_t)_{{h}/{4}}^{\tau}$
via proper finite difference formulas.  If we denote the spatial order of convergence of
$v,\,v_{t}$ as $f_{h}^{\tau}$ and $g_{h}^{\tau},$ respectively, then it follows from
a generalized Milne device \cite{sheng2023nonconventional,padgett2018quenching,liu2023semi}
that
$$
f_{h}^{\tau}(t)\approx \frac{1}{\ln2}\ln\frac{\|v_{h}^{\tau}-v_{{h}/{2}}^{\tau}\|_{2}}{\|v_{{h}/{2}}^{\tau}-v_{{h}/{4}}^{\tau}\|_{2}},~~
g_{h}^{\tau}(t)\approx \frac{1}{\ln2}\ln\frac{\|(v_{t})_{h}^{\tau}-(v_{t})_{{h}/{2}}^{\tau}\|_{2}}{\|(v_{t})_{{h}/{2}}^{\tau}-(v_{t})_{{h}/{4}}^{\tau}\|_{2}}.
$$

\begin{table}
\begin{center}
 \tabcolsep 0.06in\small
\caption{Computed orders of convergence $f_{h}^{\tau}(t), g_{h}^{\tau}(t)$ corresponding to
various $b$ ($a=2,\sigma=1.8$).}

\vspace{2mm}
\begin{tabular}{llllllllllll}
 \hline
  $b$ & $-1.0$ &$-0.9$ &$-0.8$ &$-0.7$ &$-0.6$ &$-0.5$ & $-0.4$ &$-0.3$ & $-0.2$& $-0.1$  \\\hline
  $f_{h}^{\tau}(t)$ & 1.1019 & 1.0866 & 1.0499 & 1.0420 & 1.0450 & 1.0690 & 1.0603 & 1.0515 & 1.0368 & 1.0118\\
  $g_{h}^{\tau}(t)$  &2.1836 &2.1477 &1.8637 &1.8268 &1.8891 &2.1809 &2.1537 &2.0354 &1.8466 &1.5455 \\ [2pt]\hline
 \end{tabular}
 \label{tab10}
 \end{center}
\end{table}

\begin{figure}
\begin{center}
\epsfig{file=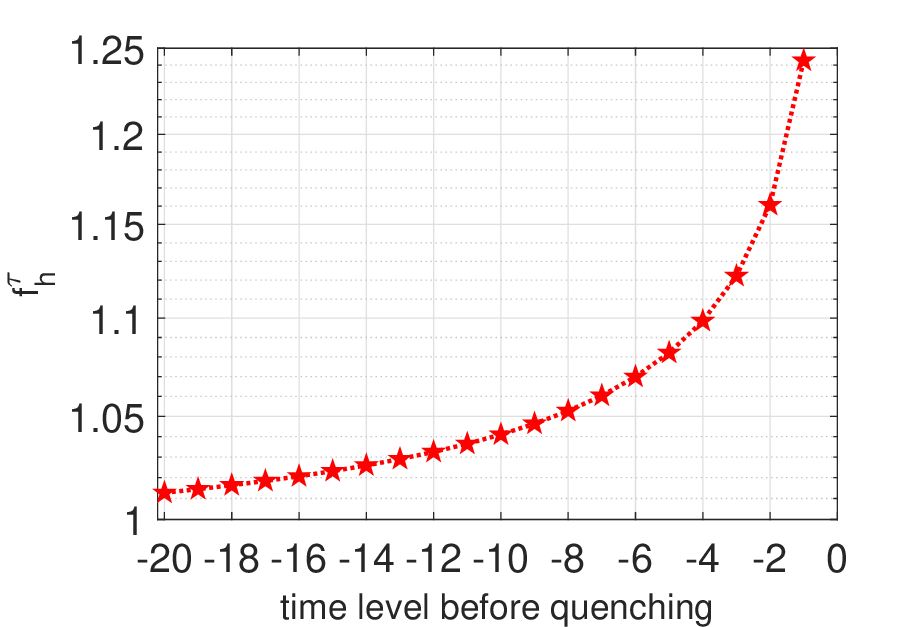,width=7.3cm,height=5.3cm}
\epsfig{file=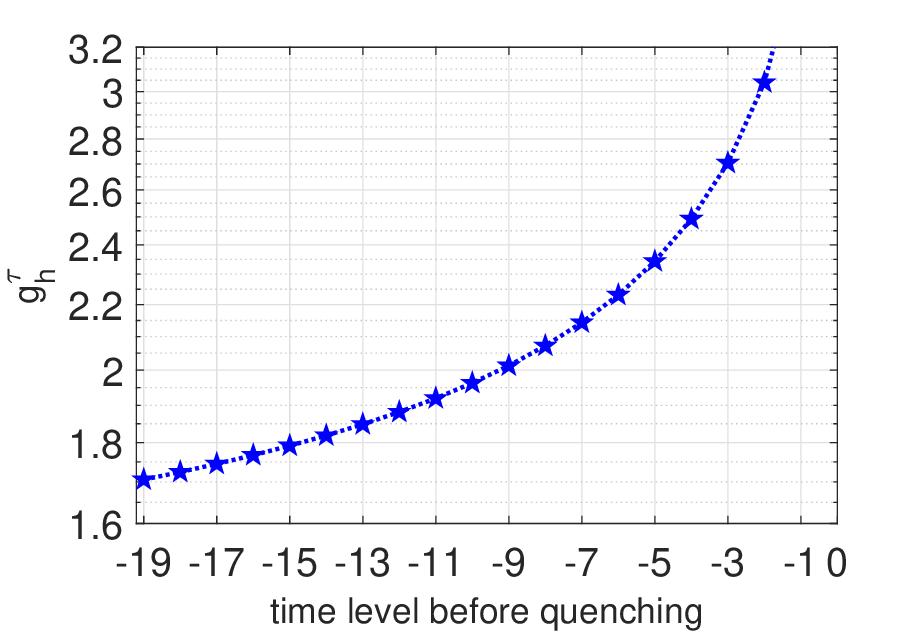,width=7.3cm,height=5.3cm}
\caption{Computed orders of convergence of $v$ (Left), $v_{t}$ (Right) within the final 20
temporal steps immediately before quenching ($a=2,\sigma=1.8,b=-0.4$).}
\label{fig11}
\end{center}
\end{figure}

Table \ref{tab10} provides estimates of the orders of convergences $f_{h}^{\tau}(t)$ and $g_{h}^{\tau}(t)$
corresponding to different quantities of $b.$ It is evident that the convection coefficient affects the convergence
slightly. As we may see, as $b$ increases, the orders of convergence $f_{h}^{\tau}(t)$ and $g_{h}^{\tau}(t)$
decreases first and then increases. Particularly speaking, when $b = -0.5,$ we have
$f_{h}^{\tau}(t)\approx 1.0690$ and $g_{h}^{\tau}(t)\approx 2.1809,$ respectively. These are consistent
with the finite difference formula (\ref{e2-9}).

At the same time, Fig. \ref{fig11} demonstrates orders of convergence of the numerical solution $v$ and
its temporal derivative $v_{t}$ when $a=2,b=-0.4.$ A fractional order, $\sigma=1.8,$ is adopted.
The simulation data are taken from the last twenty temporal steps prior to quenching.
The order of convergence of the numerical solution $v$ is found to be 1.0603 approximately. This implies
a linear convergence in both situations. The results indicate a similar order of convergence as compared
that for solving one-sided fractional order problems
\cite{padgett2018quenching,martinez2022nonlinear,liu2023semi} where the order is 1.0772. Further, the
order of convergence of $v_{t}$ is again higher at $g_{h}^{\tau}\approx 2.1537.$
The experiments may suggest that the stable semi-adaptive algorithm (\ref{e2-9}) is accurate and
highly practical for simulating quenching solutions of the singular
two-sided fractional convection-diffusion problems and their key ingredients.

\section{Concluding remarks and continuing endeavors }
\label{sec: sec5}

In this paper, we develop a semi-adaptive finite difference method for the one-dimensional, two-sided fractional order
convection-diffusion quenching problem. We precisely prove the positivity, monotonicity, and stability of the numerical
solutions. Multiple simulations illustrate the influence of various factors, such as the fractional order derivative, $b$ in
the convection term, and $\theta$ in the source term, on key parameters in quenching problems: position, time, and
critical length, using comprehensive tables and figures. Notably, when $b = 0,$ the trends of these three key
parameters are well consistent with results published in the existing literature. Additionally, we present the convergence
order of the numerical method using generalized Milne devices and examine the influence of $b$ in the convection
term on it. It is demonstrated that the proposed algorithm is highly efficient and practical for this fractional order
quenching problem.

Based on the study presented in this paper, our ongoing research will address multi-dimensional dual fractional order
quenching problems. We plan to initially focus on the influence of convection on the quenching phenomena.
Dimensional and operator splitting techniques, including the Glowinski-Le Tallec formula \cite{sk2001,lin2007finite,Glowinski},
will be employed to manage higher density computations. \color{red} Valuable physical conservations and preservations, such as free
energy preservation \cite{padgett2018quenching,martinez2022nonlinear}, will be introduced and investigated. These may
lead to solid analysis about the numerical convergence. \color{black}
More challenging multi-dimensional FPDE problems will be explored and analyzed. Physics-informed
realistic analysis and mathematical machine learning will also be carefully conducted for potential applications.

\section*{Acknowledgments}

The authors would like to express their gratitudes to the editor and anonymous reviewers for their time and extremely
valuable comments. The suggestions made by the reviewers greatly enhanced the quality and presentation of this paper.
The first and second authors are partially supported by the National Science Foundation of China (Grant No. 12062024)
and Funds for Ningxia University Scientific Research Projects (Grant No. NYG2024039).
The third author is supported in part by the National Science Foundation (Grant No. DMS-2318032) and
Simons Foundation (Grant No. MPS-1001466), USA.


\small


\bibliography{bibfile8}

\end{document}